\documentclass[12pt]{amsart}
\usepackage[english]{babel}
\usepackage{amsmath}
\usepackage{amsthm}
\usepackage{amssymb}
\usepackage{mathrsfs}
\usepackage{enumerate}
\usepackage[notcite, final, notref]{showkeys}
\usepackage{dsfont}
\usepackage{tikz}
\usepackage[siunitx]{circuitikz}
\usetikzlibrary{arrows,calc,chains,shapes,dsp}
\newtheorem{theorem}{Theorem}[section]

\newtheorem{lemma}[theorem]{Lemma}
\newtheorem{proposition}[theorem]{Proposition}
\newtheorem{observation}[theorem]{Observation}
\newtheorem{corollary}[theorem]{Corollary}

\theoremstyle{definition}
\newtheorem{remark}[theorem]{Remark}
\newtheorem{example}[theorem]{Example}
\usepackage{xcolor}

\title[Composition of rational functions and applications]
{Composition of rational functions:\\
state-space realization and applications}

\author[D. Alpay]{Daniel Alpay}
\address{(DA) 
Faculty of Mathematics, Physics, and Computation\\
Schmid College of Science and Technology\\
Chapman University\\
One University Drive
Orange, California 92866\\
USA}
\email{alpay@chapman.edu}
\thanks{Daniel Alpay thanks the Foster G. and Mary McGaw Professorship in
Mathematical Sciences, which supported this research.}

\author[I. Lewkowicz]{Izchak Lewkowicz}
\address{(IL) Department of electrical engineering 
Ben-Gurion University of the Negev\\ P.O.B. 653\\ Beer-Sheva, 84105\\
Israel}
\email{izchak@bgu.ac.il}

\begin{document}
\bibliographystyle{plain}
\begin{abstract}
We define two versions of compositions of matrix-valued rational
functions of appropriate sizes and whenever analytic at infinity,
offer a set of formulas for the corresponding state-space
realization, in terms of the realizations of the original functions.
Focusing on
positive real functions, the first composition is applied to
electrical circuits theory along with introducing a connection to
networks of feedback loops. 
The second composition is applied to
Stieltjes functions.
\end{abstract}
\maketitle


\date{today}
\tableofcontents
\section{Introduction}
\setcounter{equation}{0}

This work is focused on composition $F_L(F_R)$ of rational functions
$F_L(z)$, $F_R(z)$, where the subscript stands for ``left" and ``right".
In general
composition of functions is classical, e.g. \cite{CowMacC1995}.
Although composition of rational functions plays an important role in
the theory of dynamical systems (see e.g.
\cite{AlapayJorgLewkMarz2015}, \cite{Beard1991}), a few associated
questions are yet unsolved. We here touch upon three aspects.
\begin{itemize}

\item{}Families of functions which are closed under composition.

\item{}Applications of composition of functions.

\item{}State realization of $F_L(F_R)$ in terms of the realizations
of $F_L(z)$ and $F_R(z)$.
\end{itemize}

In principle there exist results on each of these items: For example,
applications of compositions to electrical circuits theory was studied
in \cite{Reza1984}. Little is known in the setting of realization
theory (besides the case where one of the function is a Moebius map;
see for instance \cite[Theorem 1.9, p. 35]{BarGohKaa1979}).
Realization of composition of, not necessarily rational,
\mbox{$Q$-functions} was already
addressed in \cite{Kalten1998}. 
\vskip 0.2cm

In the present work we define and study two compositions of
matrix-valued rational functions and whenever analytic at
infinity, provide formulas for the respective state-space realization.
\vskip 0.2cm

\noindent
For the first composition we offer an application to electrical
circuits. In turn, we introduce a connection with
feedback-loop networks, see the various figures below.
Implications of this novel idea go well beyond the scope of this work.
\vskip 0.2cm

\noindent
Stieltjes functions, were explored in \cite{KreinNudel1977} in the
setting of moment problems; the Nevanlinna-Pick interpolation
problem in this class was studied in that reference in the scalar 
case, and for the matrix-valued case see \cite{abgr3}, \cite{bol1,bol2}
and \cite{DyukKats1988}. We here characterize their state-space
realization and then show that both composition schemes, leave the
family of Stieltjes functions invariant. Our motivation to consider
Stieltjes functions stemmed in part from the following possible
link with statistical physics. Positive measures $\sigma$ on
$(0,\infty)$ play an important role in statistical physics, as
functions of repartition of energy levels of a particle (or, more
generally, of a system). The associated Laplace transform
\[
Z(\beta)=\int_0^\infty e^{-\beta e}d\sigma(e),
\]
assumed convergent in ${\rm Re}\, \beta>0$, is called the partition
function. See for instance \cite[p. 138]{zbMATH06649453},
\cite[p. 66]{zbMATH05917389}.
When $\sigma$ is discrete, with unit jumps at $E_1,E_2,\ldots$ we have
\begin{equation}
Z(\beta)=\sum_{j}e^{-\beta E_j}.
\label{dis}
\end{equation}
One can associate with such a measure another object, namely the
function $\varphi$ defined by
\begin{equation}
\varphi_\sigma(z)=\int_0^\infty\frac{z}{t-iz}d\sigma(t),
\end{equation}
provided $\sigma$ satisfies
\begin{equation}
\label{lll}
\int_0^\infty \frac{d\sigma(t)}{1+t}<\infty
\end{equation}
The function $\varphi_\sigma$ is a Stieltjes function and the study
of compositions of such functions, associating to two measures
$\sigma_1$ and $\sigma_2$ on $(0,\infty)$ a third measure
corresponding to the composition $\varphi_{\sigma_1}(\varphi_{\sigma_2})$
(together with possibly an imaginary constant; see formula
\eqref{eq:integral}) should have some physical interpretation, in
particular in the case of discrete finite measures.

\section{Realization of rational functions analytic at infinity}
\setcounter{equation}{0}

We first recall that a $p\times m$-valued function $F(z)$, analytic at
infinity, can be written in the form
\begin{equation}\label{eq:BasicRealization}
F(z)=D+C(zI_n-A)^{-1}B,
\end{equation}
where $D=F(\infty)$ and where $A,B,C$ are matrices of appropriate sizes.
Expression \eqref{eq:BasicRealization} is called a realization. 
Sometimes we shall find it convenient to use, the same $A,B,C,D$, the
engineering shorthand notation, (introduced by H.H. Rosenbrock,
see e.g. \cite[Chapter 1, Section 2]{Rosenbr1974}) of a $(n+p)\times(n+m)$
realization array\begin{footnote}{$R$ stands for ``realization"
or ``Rosenbrock".}
\end{footnote} $R_F$,
\begin{equation}\label{eq:RealizationArray}
R_F=
\begin{footnotesize}
\left(\begin{array}{c|c}
A&B\\
\hline
C&D\end{array}\right)
\end{footnotesize}
\end{equation}
Whenever for a given $F$, analytic at infinity, the dimension of $A$ in
\eqref{eq:BasicRealization}, \eqref{eq:RealizationArray} is the smallest
possible, the realization is called ~{\em minimal}~ and the dimension
of $A$ is ~{\em the McMillan degree of} $F(z)$. In this case, the
realization is unique up to a change of coordinates
meaning that for a $n\times n$ non-singular matrix $S$, 
\begin{equation}\label{eq:Coorninates}
\left(\begin{smallmatrix}S&0\\0&I_p\end{smallmatrix}\right)^{-1}
\underbrace{
\left(\begin{smallmatrix}A&B\\C&D\end{smallmatrix}\right)
}_{R_F}
\left(\begin{smallmatrix}S&0\\0&I_m\end{smallmatrix}\right),
\end{equation}
is another minimal realization of $F(z)$ similar to $R_F$. In particular,
the spectrum of the $A$ part, is preserved.
\vskip 0.2cm

\noindent
Up to this point, the above realization description is a classical
textbook material and we refer the reader to
\cite[Section 3.4]{AnderVongpa1973},
\cite{BarGohKaa1979},
\cite[Section 6.4 and Remark 6.7.4]{Sontag1998}.
\vskip 0.2cm

\noindent
For future reference we cite additional known results

\begin{proposition}\label{Pn:RealizationOfInverse}
Let
\[
F(z)=C(zI-A)^{-1}B+D
\]
be a $p\times p$-valued rational function, where $D$
is non-singular. Then, $\left(F(z)\right)^{-1}$ is well defined
and a realization array associated with it, i.e. with
\begin{equation}\label{eq:DefInvF}
\left(F(z)\right)^{-1}=\left(C(sI-A)^{-1}B+D\right)^{-1}
:=C_{\rm inv}(sI-A_{\rm inv})^{-1}B_{\rm inv}+D_{\rm inv}~,
\end{equation}
can be written as
\begin{equation}\label{eq:RealizationInverse}
\begin{footnotesize}
\left(\begin{array}{r|r}
A_{\rm inv}&B_{\rm inv}\\
\hline
C_{\rm inv}&D_{\rm inv}
\end{array}\right)
\end{footnotesize}
=
\begin{footnotesize}
\left(\begin{array}{c|c}
A-BD^{-1}C&-BD^{-1}\\
\hline
D^{-1}C&D^{-1}
\end{array}\right)
\end{footnotesize}.
\end{equation}
\end{proposition}
\vskip 0.2cm

For proof see e.g. \cite[Theorem 2.4]{BarGohKaaRan2010}.
\vskip 0.2cm

\begin{remark}\label{Rm:RealizeInverse}
{\rm 
One can re-write Eq. \eqref{eq:RealizationInverse} as
\[
\left(\begin{smallmatrix}
A_{\rm inv}&B_{\rm inv}\\C_{\rm inv}&D_{\rm inv}
\end{smallmatrix}\right)
=
\left(\begin{smallmatrix}
A&0\\0&0
\end{smallmatrix}\right)
+
\left(\begin{smallmatrix}
-B\\~~I_p
\end{smallmatrix}\right)
{\scriptstyle D}^{-1}
\left(\begin{smallmatrix}
C&&I_p
\end{smallmatrix}\right).
\]
}
\end{remark}
\vskip 0.2cm

\noindent
The following known result, see e.g. \cite[Section 2.5]{BarGohKaaRan2010},
will also be useful.

\begin{proposition}\label{Pn:RealizationOfProduct}
Let
\[
F_1(z)=C_1(sI-A_1)^{-1}B_1+D_1
\quad{\rm and}\quad
F_2(z)=C_2(sI-A_2)^{-1}B_2+D_2
\]
be $p\times k$ and $k\times  m$-valued rational functions, respectively.
Then, $F_1F_2$ is a $p\times m$-valued rational function whose state-space
realization may be given by
\begin{equation}\label{eq:RealizationOfProduct}
R_{F_1F_2}=
\begin{footnotesize}
\left(\begin{array}{cc|c}
A_1&B_1C_2&B_1D_2\\
0&A_2&B_2\\
\hline
C_1&D_1C_2&D_1D_2
\end{array}\right)
\end{footnotesize}.
\end{equation}
\end{proposition}
\vskip 0.2cm

\vskip 0.2cm

\noindent
We next recall in the tensor (a.k.a. Kronecker) product $M\otimes N$
of a pair of matrices $M\in\mathbb{C}^{p \times l}$ and
$N\in\mathbb{C}^{m \times q}$ so that $M\otimes N$ is of dimensions
$mp\times lq$ and takes the form\begin{footnote}
{
$M=\left(\begin{smallmatrix}
m_{11}&m_{12}&\cdots&m_{1,l}\\
m_{21}&m_{22}&\cdots&m_{2,l}\\
\vdots&\vdots&\cdots&\vdots\\
m_{p1}&m_{p2}&\cdots&m_{p,l}
\end{smallmatrix}\right)
$}
\end{footnote}

\[
M\otimes N:=\left(\begin{smallmatrix}
m_{11}N&m_{12}N&\cdots&m_{1,l}N\\
m_{21}N&m_{22}N&\cdots&m_{2,l}N\\
\vdots&\vdots&\cdots&\vdots\\
m_{p1}N&m_{p2}N&\cdots&m_{p,l}N
\end{smallmatrix}\right).
\]
For more information, see e.g.  \cite[Chapter 4]{HornJohnson2}.
\vskip 0.2cm

\noindent
We next formulate the focal problem of this work.
\vskip 0.2cm

\noindent
{\bf Problem formulation}

Let $F_L(z)$ be\begin{footnote}{Recall, the subscript stands for
``Left" and ``Right".}\end{footnote} a $p\times p$-valued rational
function of McMillan degree $~n~$ and let $F_R(z)$ be a
$q\times q$-valued rational function of McMillan degree $~m$.
Their minimal realization arrays are $(n+p)\times(n+p)$
and $(m+q)\times(m+q)$, respectively
\begin{equation}\label{eq:OrigRealiz}
R_L=
\begin{footnotesize}
\left(\begin{array}{r|r}
A_L&B_L\\
\hline
C_L&D_L
\end{array}\right)\end{footnotesize}
\quad\quad\quad
R_R=
\begin{footnotesize}
\left(\begin{array}{r|r}
A_R&B_R\\
\hline
C_R&D_R
\end{array}\right)\end{footnotesize},
\end{equation}
i.e.
\begin{equation}\label{eq:LeftRightFunctions}
F_L(z)=D_L+C_L\left(zI_n-A_L\right)^{-1}B_L~,
\quad\quad\quad
F_R(z)=D_R+C_R(zI_m-A_R)^{-1}B_R~.
\end{equation}
Assuming that $~F_L(F_R)~$ a composition of these functions, is
well defined, we seek a formula for a state-space realization
of this composition\begin{footnote}{The subscript stands for
``composition".}\end{footnote},
\begin{equation}\label{eq:FrameworkRealization}
F_L(F_R(z))=D_{\rm comp}+C_{\rm comp}(zI_k-A_{\rm comp})^{-1}B_{\rm comp}~,
\end{equation}
in terms of the realization of the original systems \eqref{eq:OrigRealiz},
for some natural $k$, i.e. a corresponding realization array is,
\begin{equation}\label{eq:DefCompFunct}
R_{F_L(F_R)}=
\begin{footnotesize}
\left(\begin{array}{r|r}
\overbrace{A_{\rm comp}}^k&B_{\rm comp}\\
\hline
C_{\rm comp}&D_{\rm comp}
\end{array}\right).\end{footnotesize}
\end{equation}
In particular, find both: $k$, see \eqref{eq:FrameworkRealization},
\eqref{eq:DefCompFunct}, the
dimension of the realization of the composed
system and the corresponding McMillan degree.

\begin{remark}{\rm
a. Even when one starts with minimal realizations of $F_L~$ and $~F_R$,
of McMillan degrees $n$ and $m$ respectively, $k$ the dimension of
the realization, see \eqref{eq:FrameworkRealization},
\eqref{eq:DefCompFunct} is not necessarily minimal,
i.e. $k$ may be bigger than
the McMillan degree of the composition.
\vskip 0.2cm

\noindent
b. In Section \ref{Sec:FirstComposition} one obtains that in
the realization, see 
\eqref{eq:FrameworkRealization}, \eqref{eq:DefCompFunct}
\[
k=mn,
\]
and in fact, this is the McMillan degree of the composed function,
see Remark \ref{Rm:degree}.\\
In contrast, in Section \ref{Sec:SecondComposition}
\[
k\leq n.
\]
}
\end{remark}

\section{Composition of functions - first version}
\label{Sec:FirstComposition}
\setcounter{equation}{0}

Presentation of the first version of composition of functions, is
split into three subcases.

\subsection{The case where $~f_R(z)~$ is scalar}

\begin{proposition}\label{Pn:f_rScalar}
Let $~f_R(z)~$ be a scalar-valued\begin{footnote}{To ease reading,
scalar functions are denoted by small letters.}\end{footnote} rational
function, of the form,
\[
f_R(z)=d_R+c_R\left(zI_m-A_R\right)^{-1}b_R~,
\]
then
\[
F_l\left(f_R(z)\right)=D_L+C_L\left(f_R(z)I_n-A_L\right)^{-1}B_L
\]
admits a state space realization of the form
\[
F_l\left(f_R(z)\right)=
D_{\rm comp}+C_{\rm comp}(zI_k-A_{\rm comp})^{-1}B_{\rm comp}~,
\]
of state dimension
\[
k=mn,
\]
and a realization array of the form
\[
\begin{footnotesize}
\left(\begin{array}{c|c}
A_{\rm comp}&B_{\rm comp}\\
\hline
C_{\rm comp}&D_{\rm comp}
\end{array}\right)\end{footnotesize}
=
\begin{footnotesize}
\left(\begin{array}{c|c}
I_n\otimes{A_R}-(I_n\otimes{b_R})
\left(d_RI_n-A_L\right)^{-1}(I_n\otimes{c_R})&
-({I_n}\otimes{b_R})\left(d_RI_n-A_L\right)^{-1}B_L\\
\hline
C_L\left(d_RI_n-A_L\right)^{-1}(I_n\otimes{c_R})
& D_L+ C_L\left(d_RI_n-A_L\right)^{-1}B_L
\end{array}\right).\end{footnotesize}
\]
\end{proposition}

{\bf Proof :}\quad
By construction,
\[
\begin{matrix}
F_L(f_R)&=&
D_L+C_L\left(f_R(z)I_n-A_L\right)^{-1}B_L
\\~\\~&=&
D_L+C_L\left(
\left(\begin{smallmatrix}
f_R(z)&~&~\\
~     &\ddots&~\\
~     &~     &f_R(z)
\end{smallmatrix}\right)-A_L\right)^{-1}B_L
\\~\\~&=&
D_L+C_L\left(
\left(\begin{smallmatrix}
d_R+c_R\left(zI_m-A_R\right)^{-1}b_R&~&~\\
~     &\ddots&~\\
~     &~     &d_R+c_R\left(zI_m-A_R\right)^{-1}b_R
\end{smallmatrix}\right)-A_L\right)^{-1}B_L
\\~\\
~&=&
D_L+C_L\left(
d_RI_n-A_L+
\left(\begin{smallmatrix}
c_R\left(zI_m-A_R\right)^{-1}b_R&~&~\\
~     &\ddots&~\\
~     &~     &c_R\left(zI_m-A_R\right)^{-1}b_R
\end{smallmatrix}\right)\right)^{-1}B_L
\\~\\
~&=&
D_L+C_L\left(
(d_RI_n-A_L)+I_n\otimes\left(c_R\left(zI_m-A_R\right)^{-1}b_R\right)
\right)^{-1}B_L
\\~\\~&=&
D_L+C_L\left(d_RI_n-A_L+(I_n\otimes c_R)\left(
I_n\otimes(zI_m-A_R)^{-1}\right)(I_n\otimes b_R)
\right)^{-1}B_L
\\~\\~&=&
D_L+C_L\left(
\underbrace{d_RI_n-A_L}_{\hat{D}}+
\underbrace{(I_n\otimes c_R)}_{\hat{C}}
\left(zI_{nm}-
\underbrace{I_n\otimes A_R}_{\hat{A}}\right)^{-1}
\underbrace{(I_n\otimes b_R)}_{\hat{B}}
\right)^{-1}B_L
\\
~&=&
D_L+C_L\left(\hat{D}^{-1}+\hat{D}^{-1}\hat{C}\left(zI_{nm}-
(\hat{A}+\hat{B}\hat{D}^{-1}\hat{C})\right)^{-1}(-\hat{B}\hat{D}^{-1})\right)B_L
\end{matrix}
\]
\[
\begin{matrix}
~&=&
\underbrace{D_L+C_L\hat{D}^{-1}B_L}_{D_{\rm comp}}+
\underbrace{C_L\hat{D}^{-1}\hat{C}}_{C_{\rm comp}}\left(zI_{nm}-
\underbrace{(\hat{A}+\hat{B}\hat{D}^{-1}\hat{C})}_{A_{\rm comp}}
\right)^{-1}\underbrace{(-\hat{B}\hat{D}^{-1}B_L)}_{B_{\rm comp}}~,
\end{matrix}
\]
where we have used Proposition \ref{Pn:RealizationOfInverse} with
\[
\hat{A}:=I_n\otimes A_R\quad\quad
\hat{B}=I_n\otimes b_R\quad\quad
\hat{C}:=I_n\otimes c_R\quad\quad
\hat{D}:=d_RI_n-A_L~,
\]
and thus the construction is complete.
\qed

\begin{remark}\label{Rm:RealizeFirstComp}
{\rm
One can re-write the last result as,
\[
\begin{footnotesize}
\left(\begin{array}{r|r}
A_{\rm comp}&B_{\rm comp}\\
\hline
C_{\rm comp}&D_{\rm comp}
\end{array}\right)
\end{footnotesize}
=
\left(\begin{smallmatrix}
I_n\otimes A_R &&0\\~\\0&&D_L\end{smallmatrix}\right)
+
\left(\begin{smallmatrix}-I_n\otimes b_R\\~\\C_L
\end{smallmatrix}\right)
\left(\begin{smallmatrix}d_RI_n-A_L\end{smallmatrix}
\right)^{-1}
\left(\begin{smallmatrix}I_n\otimes c_R&&B_L
\end{smallmatrix}\right).
\]
}
\end{remark}
\subsection{The case where $A_L$ is diagonalizable}

Here diagonalizability assumption of $A_L$, the state matrix
associated with $F_L(z)$ essentially reduces the problem to
a composition by a sum of {\em degree one} rational functions.
The details are as follows.
\vskip 0.2cm

\noindent
We start by diverting a little, and exploit diagonalizability of $A$ to
obtain a result whose applicability is well beyond the scope of this work.

\begin{proposition}\label{Pn:ProjectionRealization}
Let $A\in{\mathbb C}^{n\times n}$ be a diagonalizable matrix and
let $~a_1,~\ldots,~ a_{\nu}\in{\mathbb C}$, for some
$\nu\in[1,~n]$, be its distinct eigenvalues. Denote by
$n_1,~\ldots~,~n_{\nu}$ the corresponding algebraic multiplicity,
i.e. for some (non-unique) non-singular $V\in{\mathbb C}^{n\times n}$,
\begin{equation}\label{eq:diagonalization}
A=V^{-1}\left(\begin{smallmatrix}
a_1I_{n_1}&~         &~     & ~               \\
~         &a_2I_{n_2}& ~    & ~               \\
~         &~         &\ddots& ~               \\
~         &~         &~     &a_{\nu}I_{n_{\nu}}
\end{smallmatrix}\right)V
\quad\quad\quad\quad{\scriptstyle\sum
\limits_{j=1}^{\nu}n_j}={\scriptstyle n}~.
\end{equation}
$(I)$\quad Let $B$ a ${n\times m}$ matrix.
Then, the pair $A, B$ is controllable if and only if, with $V$
from \eqref{eq:diagonalization} one can write
\begin{equation}\label{eq:PartitioningB}
V^{-1}B=\left(\begin{smallmatrix}
{\beta_1}
\\~\\
{\beta_2}
\\
\vdots
\\~\\
{\beta_{\nu}}
\end{smallmatrix}\right)
\begin{array}{l}
\left.
\vphantom
{\begin{smallmatrix}{\scriptstyle\beta}_1
\end{smallmatrix}}\right\}{\scriptstyle n_1}
\\ 
\left.
\vphantom
{\begin{smallmatrix}{\scriptstyle\beta}_2
\end{smallmatrix}}\right\}{\scriptstyle n_2}
\\ 
\left. 
\vphantom{\begin{smallmatrix} b_o \\
b_o\end{smallmatrix}}\right.
\\
\left.
\vphantom
{\begin{smallmatrix}{\scriptstyle\beta}_{\nu}
\end{smallmatrix}}\right\}{\scriptstyle n_{\nu}}
\end{array}
\end{equation}
where each of the matrices ${\beta_1}$, $\ldots$, ${\beta_{\nu}}$ is of 
a full rank.
\vskip 0.2cm

\noindent
In particular, $m\geq\max(n_1,~n_2,~\ldots,~n_{\nu})$.
\vskip 0.2cm

\noindent
$(II)$\quad Let $C$ be a ${p\times n}$ matrix.
Then, the pair $A, C$ is observable if and only if, with
$V$ from \eqref{eq:diagonalization} one can write
\begin{equation}\label{eq:PartitioningC}
CV=(~\underbrace{{\gamma}_1}_{n_1}\quad
\underbrace{{\gamma}_2}_{n_2}\quad\cdots\quad
\underbrace{{\gamma}_{\nu}}_{n_{\nu}}~)
\end{equation}
where each of the matrices ${\gamma}_1$, $\cdots$, ${\gamma}_{\nu}$
is of a full rank.
\vskip 0.2cm

\noindent
In particular, $p\geq\max(n_1,~n_2,~\ldots,~n_{\nu})$.
\vskip 0.2cm

\noindent
$(III)$\quad Let
\[
F(z)=C(zI_n-A)^{-1}B+D,
\]
be a $p\times m$-valued rational function where $A, B, C$ are as
above. $R_F$, an $(n+p)\times(n+m)$ realization array of $F(z)$,
\[
R_F=
\begin{footnotesize}
\left(\begin{array}{c|c}
A&B\\
\hline
C&D
\end{array}\right),
\end{footnotesize}
\]
is minimal, if and only if each of the above $2\nu$ matrices,
$~{\beta}_1~,~\ldots~,~{\beta}_{\nu}$ and
$~{\gamma}_1~,~\ldots~,~{\gamma}_{\nu}$ in
\eqref{eq:PartitioningB} and
\eqref{eq:PartitioningC} respectively, is of a full-rank.
\end{proposition}

\noindent
{\bf Proof :}\quad (I)\quad By the P-B-H eigenvector's test,
see e.g. \cite[Theorem 6.2-5]{Kailath1980}, a pair \mbox{$A$, $B$}
is ~{\em uncontrollable}, if and only if (up to relabeling the
eigenvalues of $A$) there exists \mbox{$0\not=v\in{\mathbb C}^n$}
so that
\begin{equation}\label{eq:PBH}
v^*A=a_1v^*\quad\quad{\rm and}\quad\quad v^*B=0.
\end{equation}
Using \eqref{eq:diagonalization} one can write
\[
u^*:=v^*V^{-1}\quad\quad{\rm with}\quad\quad
u=\left(\begin{smallmatrix}u_1\\0\\
\vdots\\~\\0\end{smallmatrix}\right)
\quad
0\not=u_1\in{\mathbb C}^{n_1}.
\]
Substituting in
\eqref{eq:PartitioningB}, controllability means that
\[
u_1^*{\beta}_1\not=0,
\]
and since $u_1$ is arbitrary, one may conclude that
the rank of the $n_1\times m$ matrix
${\beta}_1$ is at least $n_1$. Since similar reasoning
can be applied with $~j=2,~\ldots~.~\nu$,
this part of the claim is established.
\vskip 0.2cm

\noindent
Item (II) follows from item (I) by controllability-observability
duality.
\vskip 0.2cm

\noindent
Item (III) follows from the first two items by recalling
that a realization is minimal if and only if it is both
controllable and observable.
\qed
\vskip 0.3cm

\noindent
For a diagonalizable matrix $A$, the eigenvalues-eigenspaces
description of in \eqref{eq:diagonalization} is the best known.
However, it is inherently non-unique,
i.e. one can also write,
\[
A=\left(WV\right)^{-1}
\left(\begin{smallmatrix}
a_1I_{n_1}&~         &~     & ~               \\
~         &a_2I_{n_2}& ~    & ~               \\
~         &~         &\ddots& ~               \\
~         &~         &~     &a_{\nu}I_{n_{\nu}}
\end{smallmatrix}\right)
WV
\quad\quad\quad
W:=\left(\begin{smallmatrix}
W_1&~ &~     & ~     \\
~ &W_2& ~    & ~     \\
~ &~  &\ddots& ~     \\
~ &~  &~     &W_{\nu}
\end{smallmatrix}\right)
\quad\quad
\begin{smallmatrix}
\sum\limits_{j=1}^{\nu}n_j
=
n
\\~\\
W_j\in{\mathbb C}^{n_j\times n_j}\quad{\rm non-singular}.
\end{smallmatrix}
\]
We next introduce a ~{\em unique}~ eigenvalues-eigenspaces description 
of a diagonalizable matrix $A$, to be used in the sequel. This is an
extended version of a classical 
result, see e.g. \cite[Ch. 6, Thms. 8 \& 9]{HoffKunze1967}

\begin{lemma}\label{La:DiagonalizableMatrixProjection}
Let $A\in{\mathbb C}^{n\times n}$ be a diagonalizable matrix
and let $~{a}_1,~\ldots,~ {a}_{\nu}\in{\mathbb C}$,
(with $\nu\in[1,~n]$)
be its distinct eigenvalues,
i.e. $\nu$
is the degree of the minimal polynomial associated with $A$.\\
Then, there exist (oblique) projections\begin{footnote}
{For $j=1,~\ldots~,~\nu$ the rank of $\Pi_j$ is equal to
the algebraic multiplicity of the corresponding $a_j$.
}\end{footnote},
${\scriptstyle\Pi_1}$,
$\ldots$,
${\scriptstyle\Pi_{\nu}}$ satisfying,
\[
{\scriptstyle\Pi_j\Pi_k}=
\left\{\begin{smallmatrix}\Pi_j&&j=k\\~\\
0_n&&j\not=k\end{smallmatrix}\right.
\quad\quad
\quad\quad\sum\limits_{j=1}^{\nu}{\scriptstyle\Pi_j}=I_n~,
\]
so that one can write,
\[
A=\sum\limits_{j=1}^{\nu}{\scriptstyle\Pi_ja_j}~.
\]
Furthermore, this presentation is unique.
\vskip 0.2cm

\noindent
It now follows that the pencil associated with $A$ can
be written as,
\[
\left(I_nz-A\right)^{-1}
=\sum\limits_{j=1}^{\nu}
{\scriptstyle\Pi_j}(z-a_j)^{-1}
{\scriptstyle\Pi_j}~.
\]
\end{lemma}
\vskip 0.2cm

\noindent
Note that using \eqref{eq:diagonalization} the projections
in Lemma \ref{La:DiagonalizableMatrixProjection}
are actually given by,
\[
{\scriptstyle\Pi_1}={\scriptstyle V}^{-1}\left(\begin{smallmatrix}
I_{n_1}&~         &~     & ~               \\
~      &0\cdot I_{n_2}& ~    & ~               \\
~         &~         &\ddots& ~               \\
~         &~         &~     &0\cdot I_{n_{\nu}}
\end{smallmatrix}\right){\scriptstyle V}
\quad\cdots\quad
{\scriptstyle\Pi_{\nu}}={\scriptstyle V}^{-1}\left(\begin{smallmatrix}
0\cdot I_{n_1}&~         &~     & ~               \\
~      &0\cdot I_{n_2}& ~    & ~               \\
~         &~         &\ddots& ~               \\
~         &~         &~     &I_{n_{\nu}}
\end{smallmatrix}\right){\scriptstyle V}.
\]
We can now use Lemma
\ref{La:DiagonalizableMatrixProjection} 
to obtain a convenient state space realization of a rational
function.

\begin{lemma}\label{La:ProjectionRealization}
Let $F(z)$ be a $~p\times m$-valued rational function and assume
that the associated $n\times n$ state-matrix $A$, is diagonalizable.\\
$(I)$\quad Denote by 
$~a_1,~\ldots~,~a_n$ the eigenvalues (including multiplicity) of $A$.\\
Then, there exist rank-one (oblique) projections
$~{\scriptstyle\Pi_1}$, $\ldots$
${\scriptstyle\Pi_n}$ satisfying,
\[
{\scriptstyle\Pi_j\Pi_k}=\left\{\begin{smallmatrix}{\Pi}_j&&j=k\\~\\
0_n&&j\not=k\end{smallmatrix}\right.\quad\quad
\quad\quad\sum\limits_{j=1}^n{\scriptstyle\Pi_j}=I_n~,
\]
so that 
$F(z)$ admits a unique minimal realization of the form,
\begin{equation}\label{eq:F_L}
\begin{matrix}
F(z)&=&D+C\sum\limits_{j=1}^n{\scriptstyle\Pi_j}(z-a_j)^{-1}
{\scriptstyle\Pi_j}B
\\~\\~&=&
D+\sum\limits_{j=1}^nC_j(z-a_j)^{-1}B_j
\\~\\~&=&
D+\sum\limits_{j=1}^n(z-a_j)^{-1}C_jB_j
\\~\\~&=&
D+\sum\limits_{j=1}^nC_jB_j(z-a_j)^{-1}
\end{matrix}
\quad\quad\quad
\begin{smallmatrix}
j=1,~\ldots~,~n,&
\\~\\
B_j:={\scriptstyle\Pi_j}B&
n\times m
\\~\\
\sum\limits_{j=1}^nB_j=B&
\\~\\
C_j:=C{\scriptstyle\Pi_j}&p\times n,
\\~\\
\sum\limits_{j=1}^nC_j=C. &
\end{smallmatrix}
\end{equation}
$(II)$\quad For some $\nu\in[1, n]$, denote by
$~\hat{a}_1,~\ldots~,~\hat{a}_{\nu}$, the distinct eigenvalues of $A$.\\
Then, there exist (oblique) projections\begin{footnote}
{For $j=1,~\ldots~,~\nu$ the degree of the projection $\hat{\Pi}_j$ is
equal to the algebraic multiplicity of $\hat{a}_j$.
}\end{footnote}
$~{\scriptstyle\hat{\Pi}_1}$, $\ldots$
${\scriptstyle\hat{\Pi}_{\nu}}$ satisfying,
\[
{\scriptstyle\hat{\Pi}_j\hat{\Pi}_k}=
\left\{\begin{smallmatrix}\hat{\Pi}_j&&j=k\\~\\
0_n&&j\not=k\end{smallmatrix}\right.\quad\quad
\quad\quad\sum\limits_{j=1}^{\nu}{\scriptstyle\hat{\Pi}_j}=I_n~,
\]
so that $F(z)$ admits a minimal realization of the form,
\begin{equation}\label{eq:F_Lnu}
\begin{matrix}
F(z)&=&D+C\sum\limits_{j=1}^{\nu}{\scriptstyle\hat{\Pi}_j}(z-\hat{a}_j)^{-1}
{\scriptstyle\hat{\Pi}_j}B
\\~\\~&=&
D+\sum\limits_{j=1}^{\nu}C_j(z-\hat{a}_j)^{-1}B_j
\\~\\~&=&
D+\sum\limits_{j=1}^{\nu}(z-\hat{a}_j)^{-1}C_jB_j
\\~\\~&=&
D+\sum\limits_{j=1}^{\nu}C_jB_j(z-\hat{a}_j)^{-1}
\end{matrix}
\quad\quad\quad
\begin{smallmatrix}
j=1,~\ldots~,~{\nu},&
\\~\\
B_j:={\scriptstyle\hat{\Pi}_j}B&
n\times m
\\~\\
\sum\limits_{j=1}^{\nu}B_j=B&
\\~\\
C_j:=C{\scriptstyle\hat{\Pi}_j}&p\times n,
\\~\\
\sum\limits_{j=1}^{\nu}C_j=C. &
\end{smallmatrix}
\end{equation}
Furthermore, if $~\nu=n$ or when $~m=p=1$, this presentation is unique.
\end{lemma}

\noindent
Recalling that minimality of realization is preserved under change of
coordinates, see \eqref{eq:Coorninates}, we next exploit Lemma
\ref{La:DiagonalizableMatrixProjection} to specify a minimal realization
of a rational function to be used in the sequel.
\vskip 0.2cm

\noindent
One can now apply part (I) of Lemma \ref{La:ProjectionRealization} to
$~F_L(z)~$ and consider a composition $~F_L(F_R)$.

\begin{proposition}\label{Pn:CompositionDiagonalizable}
Consider the system in the problem formulation assuming that:\\
(i) $A_L$, the $n\times n$ state-matrix associated with $F_L(z)$
is diagonalizable,\\
(ii) The eigenvalues of $A_L$ (including multiplicity) 
denoted by $a_1,~\ldots~,~a_n\in{\mathbb C}$, 
are so that the $q\times q$ matrices
\begin{equation}\label{eq:DefDelta}
{\Delta}_j:=D_R-a_jI_q\quad\quad\quad
j=1,~\ldots~,~n,
\end{equation}
are all non-singular.
\vskip 0.2cm

\noindent
In each of the three following cases of composition $F_L(F_R)$, one obtains,
a realization of as in Eqs. \eqref{eq:FrameworkRealization},
and \eqref{eq:DefCompFunct} i.e.
\[
F_L\left(F_R(z)\right)=C_{\rm comp}\left(zI_k-
A_{\rm comp}\right)^{-1}B_{\rm comp}
+D_{\rm comp},
\]
where $k=mn$ and
\begin{equation}\label{eq:DiagonalizableRealization}
A_{\rm comp}=\left(\begin{smallmatrix}
A_R-B_R{\Delta}_1^{-1}C_R&~     &~\\
~                        &\ddots&~\\
~                        & ~    &A_R-B_R{\Delta}_n^{-1}C_R
\end{smallmatrix}\right)
\quad{\rm and}\quad
D_{\rm comp}=
\begin{smallmatrix}
D_L+\sum\limits_{j=1}^nC_j{\Delta}_j^{-1}B_j~.
\end{smallmatrix}
\end{equation}
(I)\quad If $q=n$, namely $F_R(z)$ is $n\times n$-valued,
then in \eqref{eq:FrameworkRealization} and \eqref{eq:DefCompFunct}
\begin{equation}\label{eq:CompactFormulation1}
B_{\rm comp}=-
\left(\begin{smallmatrix}
B_R{\Delta}_1^{-1}B_1\\
\vdots\\~\\
B_R{\Delta}_1^{-1}B_n
\end{smallmatrix}\right)
\quad\quad\quad
C_{\rm comp}=
\left(\begin{smallmatrix}
C_1{\Delta}_1^{-1}C_R&&
\cdots&&
C_n{\Delta}_n^{-1}C_R
\end{smallmatrix}\right).
\end{equation}
(II) If $q=p$, namely both $F_L(z)$ and $F_R(z)$ are $p\times p$-valued,
then in \eqref{eq:FrameworkRealization} and \eqref{eq:DefCompFunct},
\begin{equation}\label{eq:CompactFormulation2a}
B_{\rm comp}=-
\left(\begin{smallmatrix}
B_R{\Delta}_1^{-1}C_1B_1\\
\vdots\\~\\
B_R{\Delta}_n^{-1}C_nB_n
\end{smallmatrix}\right)
\quad\quad\quad
C_{\rm comp}=
\left(\begin{smallmatrix}
{\Delta}_1^{-1}C_R&&\cdots&&
{\Delta}_n^{-1}C_R&
\end{smallmatrix}\right),
\end{equation}
or
\begin{equation}\label{eq:CompactFormulation2b}
B_{\rm comp}=-
\left(\begin{smallmatrix}
B_R{\Delta}_1^{-1}\\
\vdots\\~\\
B_R{\Delta}_n^{-1}
\end{smallmatrix}\right)\quad\quad\quad
C_{\rm comp}=
\left(\begin{smallmatrix}
C_1B_1{\Delta}_1^{-1}C_R&&
\cdots&C_nB_n{\Delta}_n^{-1}C_R&&
\end{smallmatrix}\right).
\end{equation}
(III)\quad Assume that $p=1~$ i.e. $f_L(z)$ is scalar-valued
and $~F_R(s)~$ is $~q\times q$-valued where $~q~$ is a parameter.
Here we define $~n$ scalars
\begin{equation}\label{eq:DefGamma_j}
{\scriptstyle\eta}_j:=C_jB_j\quad\quad j=1,~\ldots~,~n,
\end{equation}
and then in \eqref{eq:FrameworkRealization} and \eqref{eq:DefCompFunct},
\begin{equation}\label{eq:CompactFormulation3}
B_{\rm comp}=-
\left(\begin{smallmatrix}
B_R{\Delta}_1^{-1}\\
\vdots\\~\\
B_R{\Delta}_n^{-1}
\end{smallmatrix}\right)
\quad
\quad
\quad
C_{\rm comp}=
\left(
\begin{smallmatrix}
{\scriptstyle\eta}_1{\Delta}_1^{-1}C_R&&
\cdots&&
{\scriptstyle\eta}_n{\Delta}_n^{-1}C_R
\end{smallmatrix}
\right).
\end{equation}
\end{proposition}

\noindent
{\bf Proof of Proposition \ref{Pn:CompositionDiagonalizable}}\quad
We here find it convenient to introduce an auxiliary function $\tilde{F}_L$
and to apply \eqref{eq:F_L} to it, i.e.
\[
\tilde{F}_L(z)=D_L+\sum\limits_{j=1}^n
{\scriptstyle\gamma}_j
\left(z-a_j\right)^{-1}{\scriptstyle\beta}_j~,
\]
where the parameters ${\scriptstyle\beta}_1~\ldots~,~{\scriptstyle\beta}_n$
and
$~{\scriptstyle\gamma}_1~,~\ldots~,~{\scriptstyle\gamma}_n$
will be defined in the sequel.
\vskip 0.2cm

\noindent
We now consider realization of composition of these
functions namely, using Eq. \eqref{eq:DefDelta},
\begin{equation}\label{eq:CompositionAuxiliaryFunction}
\begin{matrix}
\tilde{F}_L(F_R)&=&D_L+\sum\limits_{j=1}^n
{\scriptstyle\gamma}_j
\left(F_R-a_jI_q\right)^{-1}{\scriptstyle\beta}_j
\\~\\~&=&D_L+\sum\limits_{j=1}^n{\scriptstyle\gamma}_j
\underbrace{\left(C_R(zI_m-A_R)^{-1}B_R+D_R-a_jI_q\right)^{-1}
}_{\left(F_R(z)-a_jI_q\right)^{-1}}{\scriptstyle\beta}_j
\\~\\~&=&D_L+\sum\limits_{j=1}^n{\scriptstyle\gamma}_j
\left(C_R(zI_m-A_R)^{-1}B_R+{\Delta}_j\right)^{-1}
{\scriptstyle\beta}_j~.
\end{matrix}
\end{equation}
Since by assumption the $q\times q$ matrices
$~{\Delta}_1,~\ldots~,~{\Delta}_n$ are all non-singular,
using Proposition \ref{Pn:RealizationOfInverse}, one can
equivalently write Eq. \eqref{eq:CompositionAuxiliaryFunction} as
\[
\begin{matrix}
\tilde{F}_L(F_R)&=&
D_L+\sum\limits_{j=1}^n
{\scriptstyle\gamma}_j
\left(
({\Delta}_j^{-1}C_R)
\left((zI_m-(A_R-B_R{\Delta}_j^{-1}C_R)\right)^{-1}
(-B_R{\Delta}_j^{-1})
+
{\Delta}_j^{-1}
\right)
{\scriptstyle\beta}_j
\\~\\
~&=&D_L+\sum\limits_{j=1}^n
{\scriptstyle\gamma}_j
{\Delta}_j^{-1}
{\scriptstyle\beta}_j
+\sum\limits_{j=1}^n\left(
{\scriptstyle\gamma}_j
{\Delta}_j^{-1}C_R\right)
\left(zI_m-(A_R-B_R{\Delta}_j^{-1}C_R)\right)^{-1}\left(
B_R{\Delta}_j^{-1}
{\scriptstyle\beta}_j
\right)~.
\end{matrix}
\]
This can be compactly written as
\begin{equation}\label{eq:AuxiliaryCompactFormulation}
R_{\tilde{F}_L(F_R)}=
\begin{footnotesize}
\left(\begin{array}{r|r}
A_{\rm comp}&B_{\rm comp}\\
\hline
C_{\rm comp}&D_{\rm comp}
\end{array}\right)
\end{footnotesize}
=
\begin{footnotesize}
\left(\begin{array}{r|r}
\begin{smallmatrix}
A_R-B_R{\Delta}_1^{-1}C_R&~&~\\
~&\ddots&~\\
~&~&A_R-B_R{\Delta}_n^{-1}C_R
\end{smallmatrix}
&
\begin{smallmatrix}
-B_R{\Delta}_1^{-1}
{\scriptstyle\beta}_1
\\
\vdots\\~\\
-B_R{\Delta}_n^{-1}
{\scriptstyle\beta}_n
\end{smallmatrix}
\\
\hline
\begin{smallmatrix}
{\scriptstyle\gamma}_1
{\Delta}_1^{-1}C_R&~&~&
\cdots&~&~&
{\scriptstyle\gamma}_n
{\Delta}_n^{-1}C_R
\end{smallmatrix}
&
\begin{smallmatrix}
D_L+\sum\limits_{j=1}^n
{\scriptstyle\gamma}_j
{\Delta}_j^{-1}
{\scriptstyle\beta}_j
\end{smallmatrix}
\end{array}\right).
\end{footnotesize}
\end{equation}
(I) To obtain \eqref{eq:CompactFormulation1} substitute
\[
{\scriptstyle\beta}_j=B_j\quad\quad
{\scriptstyle\gamma}_j=C_j\quad\quad
j=1,~\ldots~,~n.
\]
(II)(a) To obtain \eqref{eq:CompactFormulation2a} substitute
\[
{\scriptstyle\beta}_j=C_jB_j\quad\quad
{\scriptstyle\gamma}_j\equiv 1\quad\quad
j=1,~\ldots~,~n.
\]
(II)(b) To obtain \eqref{eq:CompactFormulation2b} substitute
\[
{\scriptstyle\beta}_j\equiv 1\quad\quad
{\scriptstyle\gamma}_j=C_jB_j\quad\quad
j=1,~\ldots~,~n.
\]
(III) To obtain \eqref{eq:CompactFormulation3} substitute
\[
{\scriptstyle\beta}_j\equiv 1\quad\quad
{\scriptstyle\gamma}_j=C_jB_j={\scriptstyle\eta}_j\quad\quad
j=1,~\ldots~,~n,
\]
so the construction is complete.
\mbox{}\qed\mbox{}\\

\subsection{The case where $f_L(z)$ is scalar with $A_L$ is
non-diagonalizable}

Here, one can still obtain realization of composition of
functions. However, the technical details are not as
elegant as the diagonalizable case. For simplicity of
exposition this is illustrated through an example.

\begin{example}\label{Ex:NonDiagonalizable}
{\rm
Consider the case where
\[
f_L(z)=d_L+{\scriptstyle\frac{1}{(s+a)^2}}\quad\quad\quad a>0\quad
d_L\in{\mathbb R}\quad{\rm parameters}\begin{footnote}{In Section 
\ref{Sec:Motivating applications} we shall use the fact that 
$f_L(z)$ is positive real whenever $~d_L\geq\frac{1}{8a^2}~$.
}\end{footnote}.
\]
To see that here, $A_L$ is not diagonalizable, recall that
a minimal realization of $f_L$ may be given by,
\[
R_{f_L}=
\begin{footnotesize}
\left(\begin{array}{rr|c}
-a&1&0\\0&-a&1\\
\hline
1&0&d_L
\end{array}\right)
\end{footnotesize}.
\]
Let now $F_R(z)$ be a $q\times q$-valued rational function.
Then,
the composition $f_l(F_R)$ is given by
\[
f_l(F_R)=d_LI_q+(aI_q+F_R)^{-2}
\]
and we compute, in stages, a corresponding state-space realization.
\vskip 0.2cm

\noindent
First, a realization of $F_R+aI_q$ is trivially given by
\[
R_{F_R+aI_q}=
\begin{footnotesize}
\left(\begin{array}{c|c}
A_R&B_R\\
\hline
C_R&D_R+aI_q
\end{array}\right)
\end{footnotesize}
\]
and following Proposition  \ref{Pn:RealizationOfInverse}
one has that,
\[
R_{(F_R+I_q)^{-1}}=
\begin{footnotesize}
\left(\begin{array}{c|c}
A_R-B_R(D_R+aI_q)^{-1}C_R&-B_R(D_R+aI_q)^{-1}\\
\hline
(D_R+aI_q)^{-1}C_R&(D_R+aI_q)^{-1}
\end{array}\right).
\end{footnotesize}
\]
Next, following Proposition \ref{Pn:RealizationOfProduct}
\[
R_{(F_R+aI_q)^{-2}}=
\begin{footnotesize}
\left(\begin{array}{cc|c}
A_R-B_R(D_R+aI_q)^{-1}C_R&-B_R(D_R+aI_q)^{-2}C_R
&-B_R(D_R+aI_q)^{-2}\\
0&A_R-B_R(D_R+aI_q)^{-1}C_R&-B_R(D_R+aI_q)^{-1}\\
\hline
(D_R+aI_q)^{-1}C_R&(D_R+aI_q)^{-2}C_R&(D_R+aI_q)^{-2}
\end{array}\right)
\end{footnotesize}
\]
and finally,
\[
R_{f_L(F_R)}=
\begin{footnotesize}
\left(\begin{array}{cc|c}
A_R-B_R(D_R+aI_q)^{-1}C_R&-B_R(D_R+aI_q)^{-2}C_R&-B_R(D_R+aI_q)^{-2}\\
0&A_R-B_R(D_R+aI_q)^{-1}C_R&-B_R(D_R+aI_q)^{-1}\\
\hline
(D_R+aI_q)^{-1}C_R&(D_R+aI_q)^{-2}C_R&d_LI_q+(D_R+aI_q)^{-2}
\end{array}\right)
\end{footnotesize}.
\]
}
\qed
\end{example}
\vskip 0.2cm

\begin{remark}\label{Rm:degree}
{\rm
In Propositions \ref{Pn:f_rScalar} and \ref{Pn:CompositionDiagonalizable}
and in Example \ref{Ex:NonDiagonalizable}, the degree of the realization
of $F_L(F_R)$, the composition of functions, is (up to minimality) equal
to the product of the McMillan degrees of the original realizations.
}
\end{remark}

\section{Applications to electrical circuits and to
feedback-loop networks}
\label{Sec:Motivating applications}
\setcounter{equation}{0}

In the sequel we shall denote by\begin{footnote}
{As before, the subscript stands for ``left" or ``right"}
\end{footnote} $\mathbb{C}_L$, $\mathbb{C}_R$
the open left, right, halves of the complex plane
(and by $\overline{\mathbb C}_R$
the closed right half of the complex plane).
\vskip 0.2cm

\noindent
We shall denote by ($\overline{\mathbf P}_k$) $\mathbf{P}_k$ the set
of $k\times k$ positive (semi) definite matrices.
\vskip 0.2cm

\noindent
Recall that a $p\times p$-valued functions $F(z)$
is said to be ~{\em positive}~ if
\begin{equation}\label{eq:Positive}
\forall z\in\mathbb{C}_r\quad\quad
\begin{matrix}
F(z)\quad{\rm analytic}
\\~\\
\left(F(z)+(F(z))^*\right)
\in\mathbf{P}_p~.
\end{matrix}
\end{equation}
In engineering it is further restricted so that
\[
F(z)_{|_{z\in\mathbb{R}}}\in\mathbb{R}^{p\times p},
\]
and then called ~{\em positive real}.
For details see e.g. \cite{AnderVongpa1973}, \cite{Belev1968}, \cite{Brune1},
\cite{Brune2}, \cite{Cauer1926}, \cite{Cauer1932}, \cite{Chen1964},
\cite{Wohl1969}.
\vskip 0.2cm

We first establish a connection with the previous section.

\begin{observation}\label{Ob:PrFirstComp}
Whenever both $F_L(z)$ and $F_R(z)$ are positive real, then in each
of the three above cases, i.e. Propositions \ref{Pn:f_rScalar},
\ref{Pn:CompositionDiagonalizable} and Example
\ref{Ex:NonDiagonalizable},
the resulting composed function $F_L(F_R)$, is positive real.
\end{observation}

Indeed, in terminology of scalar functions, a positive real function
maps the right half plane to itself.
\vskip 0.2cm

\noindent
Duality between rational positive real functions and the driving point
immittance of $~R-L-C~$ electrical circuits, has already been 
recognized for about ninety years, e.g. \cite{Brune1}, \cite{Brune2},
\cite{Cauer1926}, \cite{Cauer1932}. This has lead to rich and
well-established theory, see e.g. \cite{AnderVongpa1973}, \cite{Belev1968},
\cite{Chen1964}, \cite{Wohl1969}.
\vskip 0.2cm

\begin{figure}[ht!]
  \begin{tikzpicture}[scale=1.5]
  \draw[color=black, thick]
        (0,0) to [short,o-] (2.5,0){} 
        (-0.1,1) node[]{\large{$\mathbf{Z_{\rm in}~~\rightarrow}$}}
        (0,2) to [short,o-] (2.5,2)
        (1.5,0) to [L, l=$L$,*-*] (1.5,2)
        (2.5,0) to [C, l=$C$,*-*] (2.5,2)
        ;
  \end{tikzpicture}
  \caption{${\rm Z}_{\rm in}(z)=\left((zL)^{-1}+zC\right)^{-1}.$}
  \label{Fig:CircuitFeedbackLoop1}
\end{figure}
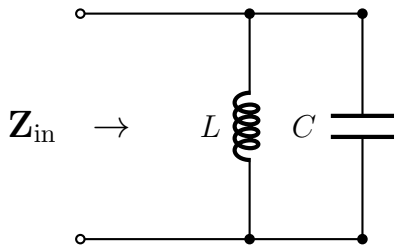

\noindent 
This duality is illustrated through two simple examples in Figures
\ref{Fig:CircuitFeedbackLoop1} and
\ref{Fig:CircuitFeedbackLoop}.
\vskip 0.2cm

\vskip 0.2cm

\begin{figure}[ht!]
\begin{tikzpicture}[scale=1.5]
\draw[color=black, thick]
       (0,0) to [short,o-] (5.3,0){} 
        (0.0,1) node[]{\large{$\mathbf{Z_{\rm in}~~\rightarrow}$}}
        (5.3,2)   to node[short]{} (1.8,2)
        (6.3,2)   to node[short]{} (10.1,2)
        (6.3,0)   to node[short]{} (10.1,0)
        (6.8,2) to [R=$R_{n-1}$,*-*] (6.8,0)
        (7.9,2) to [C=$C_{n-1}$, *-*] (7.9,0)
        (9.1,2) to [R=$R_n$,*-*] (9.1,0)
        (10.1,2) to [C=$C_n$, *-*] (10.1,0)
        (2.8,0) to [C,l=$C_1$, *-*] (2.8,2)
        (1.8,0) to [R=$R_1$, *-*] (1.8,2)
        (3.8,0) to [R=$R_2$, *-*] (3.8,2)
        (0,2) to [R=$R_L$, *-*] (1.8,2)
        (4.8,0) to [C,l=$C_2$, *-*] (4.8,2)
[color=black, thick] (5.3,0.0)node {$\bullet$}
[color=black, thick] (5.3,2.0)node {$\bullet$}
[color=black, thick] (5.55,0.0)node {$\bullet$}
[color=black, thick] (5.55,2.0)node {$\bullet$}
[color=black, thick] (5.8,0.0)node {$\bullet$}
[color=black, thick] (6.05,0.0)node {$\bullet$}
[color=black, thick] (6.3,0.0)node {$\bullet$}
[color=black, thick] (6.05,2.0)node {$\bullet$}
[color=black, thick] (6.3,2.0)node {$\bullet$}
[color=black, thick] (5.8,2.0)node {$\bullet$}
         ;
\end{tikzpicture}
\caption{\mbox{$
{\rm Z}_{\rm in}(z)={\scriptstyle R_L}+\sum\limits_{j=1}^n
\left({\scriptstyle R_j}^{-1}+z{\scriptstyle C_j}\right)^{-1}_{|_{R_L=d_L~~
R_j=\frac{{\gamma}_j}{a_j}~~C_j=\frac{1}{{\gamma}_j}}}={\scriptstyle d_L}
+\sum\limits_{j=1}^n{\scriptstyle {\gamma}_j}(z+a_j)^{-1}.$}}
\label{Fig:CircuitFeedbackLoop}
\end{figure}
\vskip 0.2cm

\noindent 
One can next address a higher level of this duality between positive real 
rational functions and the driving point immittance of $~R-L-C~$
electrical circuits: Composition of rational functions is translated, in
circuits language, to substituting elements by sub-networks, while
preserving the original configuration. For instance, $zL$ and $zC$ in
Figure \ref{Fig:CircuitFeedbackLoop1} are substituted in Figure
\ref{Fig:CircuitFeedbackLoop4} by the {\em impedance network} $Z_G$ and
the {\em admittance network} $Y_F$, respectively. This suggests
constructing an elaborate network when the basic building blocks are
$~p\times p$-valued positive real functions, see Section \ref{Sec:ConclRem}.

\begin{figure}[ht!]
  \begin{tikzpicture}[scale=1.5]
  \draw[color=black, thick]
        (0,0) to [short,o-] (2.5,0){} 
        (-0.1,1) node[]{\large{$\mathbf{Z_{\rm in}~~\rightarrow}$}}
        (0,2) to [short,o-] (2.5,2)
        (1.5,0) to [short] (1.5,0.6)
        (1.70,0.6) to [short] (1.70,1.4)
        (1.30,0.6) to [short] (1.30,1.4)
        (1.30,0.6) to [short] (1.70,0.6)
        (1.30,1.4) to [short] (1.70,1.4)
        (1.5,1.4) to [short] (1.5,2.0)
        (2.5,0) to [short] (2.5,0.6)
        (2.70,0.6) to [short] (2.70,1.4)
        (2.30,0.6) to [short] (2.30,1.4)
        (2.30,0.6) to [short] (2.70,0.6)
        (2.30,1.4) to [short] (2.70,1.4)
        (2.5,1.4) to [short] (2.5,2.0)
        ;
\draw[color=black, thick] (1.5,1)node {$\mathbf{Z_G}$};
\draw[color=black, thick] (2.5,1)node {$\mathbf{Y_F}$};
  \end{tikzpicture}
  \caption{${\rm Z}_{\rm in}(z)=\left(Y_F+{Z_G}^{-1}\right)^{-1}.$}
  \label{Fig:CircuitFeedbackLoop4}
\end{figure}
\vskip 0.2cm

\begin{figure}[ht!]
\begin{tikzpicture}
\matrix (diag_mat) [nodes={draw, fill=blue!10}, row sep=7mm, column sep=7mm]
            {
     \node[dspnodeopen,dsp/label=left] (u) {${\rm In}$};    &
     \node[dspadder,label={above left:$+$}, label={below right:$-$}] (sum) {}; &
     \node[dspfilter]    (F) {$G$};&
     ֿ\coordinate (splt) {}; &
     \node[dspnodeopen,dsp/label=above] (y) {${\rm Out}$};  \\  
         ~;&
     \coordinate (zr) ; &
     \node[dspfilter]    (H) {$F$};&
     \coordinate (zl) {}; &
                ~ ;\\
           };
     \draw[dspconn] (u) -- (sum);
     \draw[dspconn] (sum) -- (F) -- (splt) -- (y);
     \draw[dspconn] (splt) -- (zl) -- (H);
     \draw[dspconn] (H) -- (zr) -- (sum);
     \end{tikzpicture}
     \caption{~~${\rm Out}=\left(F+G^{-1}\right)^{-1}\cdot{\rm In}$}
     \label{Figure:FeedbackLoop1}
\end{figure}
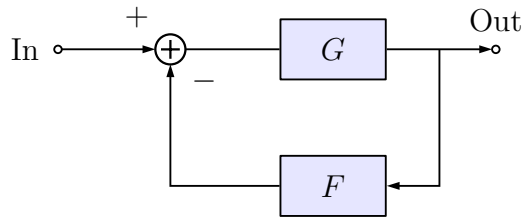

\noindent
To the above mentioned duality we now add a third aspect, namely~ {\em
interconnection of feedback loops}. This is best illustrated by an
example. Let $G$ and $F$ be square matrix-valued (possibly scalar)
rational functions so that
\[
{\rm det}(G)\not\equiv 0
\quad
\quad
{\rm and}
\quad
\quad
{\rm det}(G^{-1}+F)\not\equiv 0.
\]
Then, the input-output relation of the feedback loop in Figure
\ref{Figure:FeedbackLoop1} is well defined and is given by
\[
{\rm Out}=\left(F+G^{-1}\right)^{-1}\cdot{\rm In}.
\]

\noindent
Now, one can identify $~F~$ and $~G~$
in Figure \ref{Figure:FeedbackLoop1}, with $Y_F$ and $Z_G$
respectively, from Figure \ref{Fig:CircuitFeedbackLoop4}.
\vskip 0.2cm

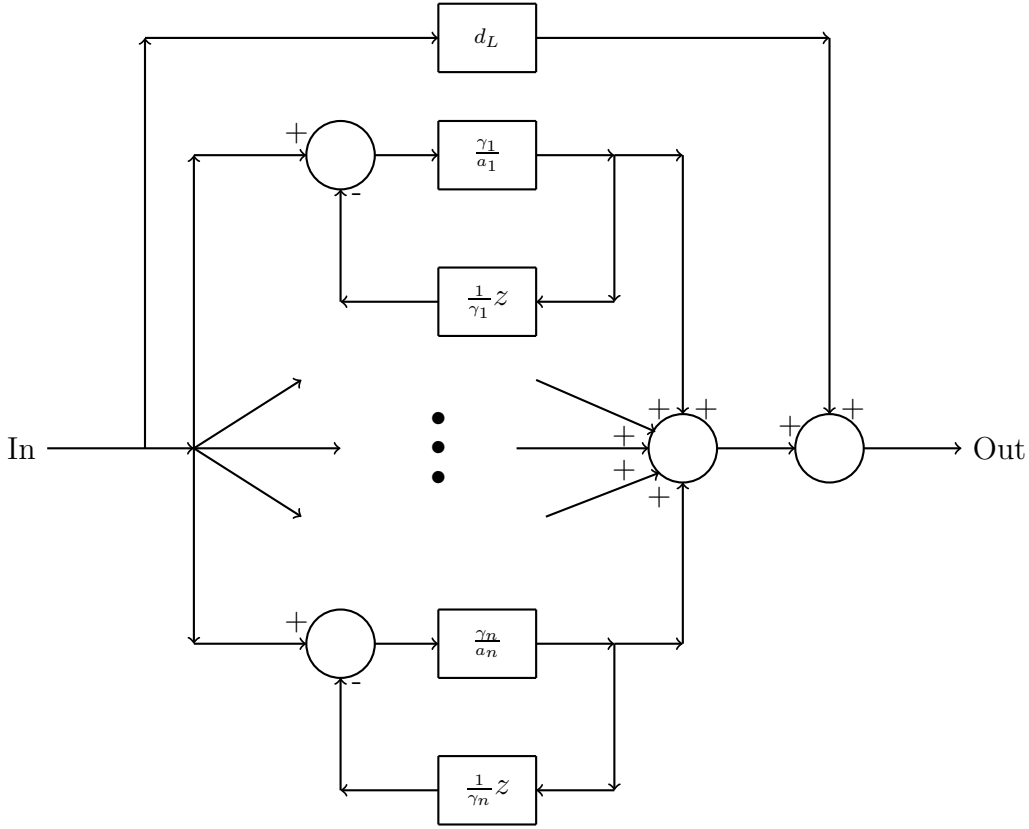
\begin{figure}[ht!]
\begin{tikzpicture}[scale=1.3]
\draw[color=black, thick] [->] (0,5)  node[left]{In} -- (1.5,5){};
\draw[color=black, thick] [->] (1.5,5) -- (1.5,8){};
\draw[color=black, thick] [->] (1.0,5) -- (1.0,9.2){};
\draw[color=black, thick] [->] (1.0,9.2) -- (4.0,9.2) {};
\draw[color=black, thick] [->] (5.0,9.2) -- (8.0,9.2) {};
\draw[color=black, thick] [->] (8.0,9.2) -- (8.0,5.35) {};
\draw[color=black, thick] [->] (1.5,5) -- (2.6,5.7){};
\draw[color=black, thick] [->] (1.5,5) -- (2.6,4.3){};
\draw[color=black, thick] [->] (5.1,4.3) -- (6.26,4.75){};
\draw[color=black, thick] [->] (4.8,5) -- (6.15,5){};
\draw[color=black, thick] [->] (5.0,5.7) -- (6.22,5.17){};
\draw[color=black, thick] [->] (1.5,5) -- (3.0,5){};
\draw[color=black, thick] [->] (1.5,5) -- (1.5,3){};
\draw[color=black, thick] [->] (1.5,3) -- (2.65,3){};
\draw[color=black, thick] [->] (1.5,8) -- (2.65,8){};
\draw[color=black, thick] [->] (3.35,8) -- (4.0,8){};
\draw[color=black, thick] [->] (3.35,3) -- (4.0,3){};
\draw[color=black, thick] [->] (5.0,8) -- (5.8,8){};
\draw[color=black, thick] [->] (5.8,8) -- (6.5,8){};
\draw[color=black, thick] [->] (5.0,3) -- (5.8,3){};
\draw[color=black, thick] [->] (5.8,3) -- (6.5,3){};
\draw[color=black, thick] [->] (6.5,3) -- (6.5,4.65){};
\draw[color=black, thick] [->] (6.85,5) -- (7.65,5){};
\draw[color=black, thick] [->] (8.35,5) -- (9.35,5)node[right] {Out};
\draw[color=black, thick] [->] (6.5,8) -- (6.5,5.35){};
\draw[color=black, thick] [->] (5.8,1.5) -- (5.0,1.5){};
\draw[color=black, thick] [->] (5.8,3.0) -- (5.8,1.5){};
\draw[color=black, thick] [->] (5.8,8.0) -- (5.8,6.5){};
\draw[color=black, thick] [->] (5.8,6.5) -- (5.0,6.5){};
\draw[color=black, thick] [->] (4.0,1.5) -- (3.0,1.5){};
\draw[color=black, thick] [->] (3.0,1.5) -- (3.0,2.65){};
\draw[color=black, thick] [->] (3.0,6.5) -- (3.0,7.65){};
\draw[color=black, thick] [->] (4.0,6.5) -- (3.0,6.5){};
\draw[color=black, thick] (3.0,3) circle (0.35){};
\draw[color=black, thick] (6.5,5) circle (0.35){};
\draw[color=black, thick] (8.0,5) circle (0.35){};
\draw[color=black, thick] (3.0,8) circle (0.35){};
\draw[color=black, thick] [-]  (4,2.65) -- (4,3.35){};
\draw[color=black, thick] [-]  (4,2.65) -- (5,2.65){};
\draw[color=black, thick] [-]  (5,3.35) -- (5,2.65){};
\draw[color=black, thick] [-]  (4,3.35) -- (5,3.35){};
\draw[color=black, thick] [-]  (4,8.35) -- (4,7.65){};
\draw[color=black, thick] [-]  (4,7.65) -- (5,7.65){};
\draw[color=black, thick] [-]  (4,8.35) -- (5,8.35){};
\draw[color=black, thick] [-]  (5,8.35) -- (5,7.65){};
\draw[color=black, thick] [-]  (4,6.85) -- (5,6.85){};
\draw[color=black, thick] [-]  (4,6.15) -- (5,6.15){};
\draw[color=black, thick] [-]  (5,1.85) -- (5,1.15){};
\draw[color=black, thick] [-]  (5,6.85) -- (5,6.15){};
\draw[color=black, thick] [-]  (4,6.85) -- (4,6.15){};
\draw[color=black, thick] [-]  (4,1.85) -- (5,1.85){};
\draw[color=black, thick] [-]  (4,1.15) -- (5,1.15){};
\draw[color=black, thick] [-]  (4,1.85) -- (4,1.15){};
\draw[color=black, thick] [-]  (4,9.55) -- (5,9.55){};
\draw[color=black, thick] [-]  (4,8.85) -- (5,8.85){};
\draw[color=black, thick] [-]  (4,8.85) -- (4,9.55){};
\draw[color=black, thick] [-]  (5,8.85) -- (5,9.55){};
\draw[color=black, thick] (4.5,3)node {${\scriptstyle\frac{{\gamma}_n}{a_n}}$};
\draw[color=black, thick] (4.5,8)node {${\scriptstyle\frac{{\gamma}_1}{a_1}}$};
\draw[color=black, thick] (4.5,1.5)node {${\scriptstyle\frac{1}{{\gamma}_n}}z$};
\draw[color=black, thick] (4.5,6.5)node {${\scriptstyle\frac{1}{{\gamma}_1}}z$};
\draw[color=black, thick] (4.5,9.2)node {${\scriptstyle d_L}$};
\draw[color=black, thick] (4.0,5.0)node {$\bullet$};
\draw[color=black, thick] (4.0,5.3)node {$\bullet$};
\draw[color=black, thick] (4.0,4.7)node {$\bullet$};
\draw[color=black, thick] (2.55,3.0)node[above] {+};
\draw[color=black, thick] (2.55,8.0)node[above] {+};
\draw[color=black, thick] (3.0,2.60)node[right] {-};
\draw[color=black, thick] (3.0,7.60)node[right] {-};
\draw[color=black, thick] (6.5,5.40)node[left] {+};
\draw[color=black, thick] (6.5,5.40)node[right] {+};
\draw[color=black, thick] (7.6,5)node[above] {+};
\draw[color=black, thick] (5.9,5)node[below] {+};
\draw[color=black, thick] (5.9,4.9)node[above] {+};
\draw[color=black, thick] (6.5,4.50)node[left] {+};
\draw[color=black, thick] (8.0,5.40)node[right] {+};
\end{tikzpicture}
\caption{$
{\rm Out}=f_L(z)\cdot
{\rm In}$\quad with\quad
$f_L(z)={\scriptstyle d_L}+\sum\limits_{j=1}^n
{\scriptstyle\gamma}_j(z+a_j)^{-1}$
}
\label{Fig:f_LFeedbackLoop}
\end{figure}

\noindent
As an engineering application of item (III) of Proposition
\ref{Pn:CompositionDiagonalizable}
see Figures \ref{Fig:CircuitFeedbackLoop},
\ref{Fig:f_LFeedbackLoop}, \ref{Fig:CompositionFeedbackLoop}.
\vskip 0.2cm

\begin{figure}[ht!]
\begin{tikzpicture}[scale=1.3]
\draw[color=black, thick] [->] (0,3)  node[left]{In} -- (0.8,3){};
\draw[color=black, thick] [->] (0.8,3) -- (0.8,4.5){};
\draw[color=black, thick] [->] (0.8,4.5) -- (4.5,4.5){};
\draw[color=black, thick] [->] (5.5,4.5) -- (9.9,4.5){};
\draw[color=black, thick] [->] (9.9,4.5) -- (9.9,3.35){};
\draw[color=black, thick] [->] (0.8,3) -- (1.65,3){};
\draw[color=black, thick] [->] (2.35,3) -- (3.0,3){};
\draw[color=black, thick] [->] (4.0,3) -- (4.8,3){};
\draw[color=black, thick] [->] (8.0,3) -- (8.8,3){};
\draw[color=black, thick] [->] (8.8,3) -- (9.55,3){};
\draw[color=black, thick] [->] (6.35,3) -- (7.0,3){};
\draw[color=black, thick] [->] (4.8,3) -- (5.65,3){};
\draw[color=black, thick] [->] (10.25,3) -- (11.00,3)node[right] {Out};
\draw[color=black, thick] [->] (4.8,1.5) -- (4.0,1.5){};
\draw[color=black, thick] [->] (8.8,1.5) -- (8.0,1.5){};
\draw[color=black, thick] [->] (4.8,3.0) -- (4.8,1.5){};
\draw[color=black, thick] [->] (8.8,3.0) -- (8.8,1.5){};
\draw[color=black, thick] [->] (3.0,1.5) -- (2.0,1.5){};
\draw[color=black, thick] [->] (7.0,1.5) -- (6.0,1.5){};
\draw[color=black, thick] [->] (2.0,1.5) -- (2.0,2.65){};
\draw[color=black, thick] [->] (6.0,1.5) -- (6.0,2.65){};
\draw[color=black, thick] (2.0,3) circle (0.35){};
\draw[color=black, thick] (6.0,3) circle (0.35){};
\draw[color=black, thick] (9.9,3) circle (0.35){};
\draw[color=black, thick] [-]  (4,2.65) -- (4,3.35){};
\draw[color=black, thick] [-]  (7,2.65) -- (7,3.35){};
\draw[color=black, thick] [-]  (8,2.65) -- (8,3.35){};
\draw[color=black, thick] [-]  (4,2.65) -- (3,2.65){};
\draw[color=black, thick] [-]  (8,2.65) -- (7,2.65){};
\draw[color=black, thick] [-]  (3,3.35) -- (3,2.65){};
\draw[color=black, thick] [-]  (4,3.35) -- (3,3.35){};
\draw[color=black, thick] [-]  (8,3.35) -- (7,3.35){};
\draw[color=black, thick] [-]  (4.5,4.85) -- (4.5,4.15){};
\draw[color=black, thick] [-]  (5.5,4.85) -- (5.5,4.15){};
\draw[color=black, thick] [-]  (5.5,4.85) -- (4.5,4.85){};
\draw[color=black, thick] [-]  (4.5,4.15) -- (5.5,4.15){};
\draw[color=black, thick] [-]  (3,1.85) -- (3,1.15){};
\draw[color=black, thick] [-]  (4,1.85) -- (3,1.85){};
\draw[color=black, thick] [-]  (8,1.85) -- (7,1.85){};
\draw[color=black, thick] [-]  (4,1.15) -- (3,1.15){};
\draw[color=black, thick] [-]  (8,1.15) -- (7,1.15){};
\draw[color=black, thick] [-]  (4,1.85) -- (4,1.15){};
\draw[color=black, thick] [-]  (8,1.85) -- (8,1.15){};
\draw[color=black, thick] [-]  (7,1.85) -- (7,1.15){};
\draw[color=black, thick] (3.5,3)node {$\scriptstyle\frac{1}{a}I_q$};
\draw[color=black, thick] (7.5,3)node {$\scriptstyle\frac{1}{a}I_q$};
\draw[color=black, thick] (3.5,1.5)node {$ F_R$};
\draw[color=black, thick] (7.5,1.5)node {$ F_R$};
\draw[color=black, thick] (5.0,4.5)node {${\scriptstyle d_L}I_q$};
\draw[color=black, thick] (1.52,3.0)node[above] {+};
\draw[color=black, thick] (2.0,2.5)node[right] {-};
\draw[color=black, thick] (6.0,2.50)node[right] {-};
\draw[color=black, thick] (5.52,3)node[above] {+};
\draw[color=black, thick] (9.43,3)node[above] {+};
\draw[color=black, thick] (9.9,3.47)node[right] {+};
\end{tikzpicture}
\caption{
${\rm Out}=f_L\left(F_R(z)\right)\cdot{\rm In}~$ with
$~f_L(z)={\scriptstyle d_L}+{\scriptstyle\frac{1}{(s+a)^2}}~~
{\scriptstyle a}>0,~{\scriptstyle d_L\in{\mathbb R}}$
parameters.}
\label{Fig:CompositionNonDiagonalizableFeedbackLoop}
\end{figure}
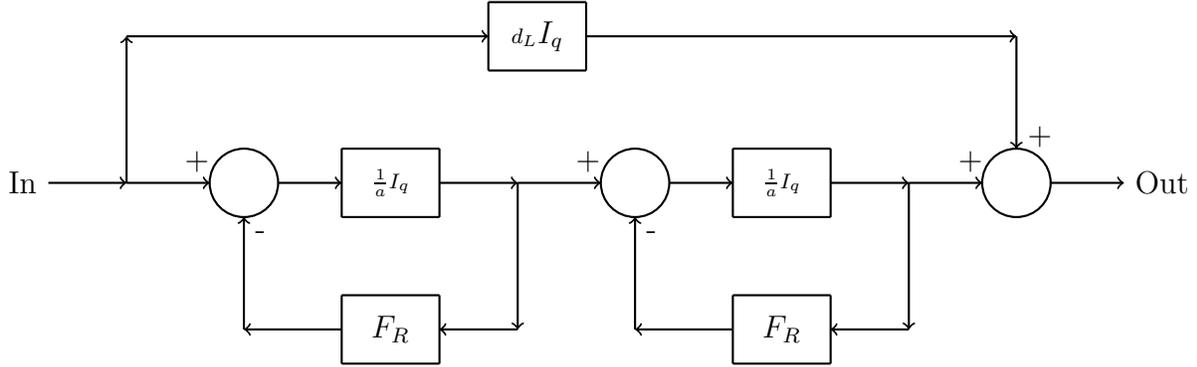

\noindent
Figure \ref{Fig:CompositionNonDiagonalizableFeedbackLoop}
presents an engineering application of Example \ref{Ex:NonDiagonalizable}.
\vskip 0.2cm

\noindent
Connection between positive real rational functions and
feedback loops is further elaborated on in Section \ref{Sec:ConclRem}.

\begin{figure}[ht!]
\begin{tikzpicture}[scale=1.3]
\draw[color=black, thick] [->] (0,5)  node[left]{In} -- (1.5,5){};
\draw[color=black, thick] [->] (1.5,5) -- (1.5,8){};
\draw[color=black, thick] [->] (1.0,5) -- (1.0,9.2){};
\draw[color=black, thick] [->] (1.0,9.2) -- (4.0,9.2) {};
\draw[color=black, thick] [->] (5.0,9.2) -- (8.0,9.2) {};
\draw[color=black, thick] [->] (8.0,9.2) -- (8.0,5.35) {};
\draw[color=black, thick] [->] (1.5,5) -- (2.6,5.7){};
\draw[color=black, thick] [->] (1.5,5) -- (2.6,4.3){};
\draw[color=black, thick] [->] (5.1,4.3) -- (6.26,4.75){};
\draw[color=black, thick] [->] (4.8,5) -- (6.15,5){};
\draw[color=black, thick] [->] (5.0,5.7) -- (6.22,5.17){};
\draw[color=black, thick] [->] (1.5,5) -- (3.0,5){};
\draw[color=black, thick] [->] (1.5,5) -- (1.5,3){};
\draw[color=black, thick] [->] (1.5,3) -- (2.65,3){};
\draw[color=black, thick] [->] (1.5,8) -- (2.65,8){};
\draw[color=black, thick] [->] (3.35,8) -- (4.0,8){};
\draw[color=black, thick] [->] (3.35,3) -- (4.0,3){};
\draw[color=black, thick] [->] (5.0,8) -- (5.8,8){};
\draw[color=black, thick] [->] (5.8,8) -- (6.5,8){};
\draw[color=black, thick] [->] (5.0,3) -- (5.8,3){};
\draw[color=black, thick] [->] (5.8,3) -- (6.5,3){};
\draw[color=black, thick] [->] (6.5,3) -- (6.5,4.65){};
\draw[color=black, thick] [->] (6.85,5) -- (7.65,5){};
\draw[color=black, thick] [->] (8.35,5) -- (9.35,5)node[right] {Out};
\draw[color=black, thick] [->] (6.5,8) -- (6.5,5.35){};
\draw[color=black, thick] [->] (5.8,1.5) -- (5.0,1.5){};
\draw[color=black, thick] [->] (5.8,3.0) -- (5.8,1.5){};
\draw[color=black, thick] [->] (5.8,8.0) -- (5.8,6.5){};
\draw[color=black, thick] [->] (5.8,6.5) -- (5.0,6.5){};
\draw[color=black, thick] [->] (4.0,1.5) -- (3.0,1.5){};
\draw[color=black, thick] [->] (3.0,1.5) -- (3.0,2.65){};
\draw[color=black, thick] [->] (3.0,6.5) -- (3.0,7.65){};
\draw[color=black, thick] [->] (4.0,6.5) -- (3.0,6.5){};
\draw[color=black, thick] (3.0,3) circle (0.35){};
\draw[color=black, thick] (6.5,5) circle (0.35){};
\draw[color=black, thick] (3.0,8) circle (0.35){};
\draw[color=black, thick] (8.0,5) circle (0.35){};
\draw[color=black, thick] [-]  (4,2.65) -- (4,3.35){};
\draw[color=black, thick] [-]  (4,2.65) -- (5,2.65){};
\draw[color=black, thick] [-]  (5,3.35) -- (5,2.65){};
\draw[color=black, thick] [-]  (4,3.35) -- (5,3.35){};
\draw[color=black, thick] [-]  (4,8.35) -- (4,7.65){};
\draw[color=black, thick] [-]  (4,7.65) -- (5,7.65){};
\draw[color=black, thick] [-]  (4,8.35) -- (5,8.35){};
\draw[color=black, thick] [-]  (5,8.35) -- (5,7.65){};
\draw[color=black, thick] [-]  (4,6.85) -- (5,6.85){};
\draw[color=black, thick] [-]  (4,6.15) -- (5,6.15){};
\draw[color=black, thick] [-]  (5,1.85) -- (5,1.15){};
\draw[color=black, thick] [-]  (5,6.85) -- (5,6.15){};
\draw[color=black, thick] [-]  (4,6.85) -- (4,6.15){};
\draw[color=black, thick] [-]  (4,1.85) -- (5,1.85){};
\draw[color=black, thick] [-]  (4,1.15) -- (5,1.15){};
\draw[color=black, thick] [-]  (4,1.85) -- (4,1.15){};
\draw[color=black, thick] [-]  (4,9.55) -- (5,9.55){};
\draw[color=black, thick] [-]  (4,8.85) -- (5,8.85){};
\draw[color=black, thick] [-]  (4,8.85) -- (4,9.55){};
\draw[color=black, thick] [-]  (5,8.85) -- (5,9.55){};
\draw[color=black, thick] (4.5,3)node {${\scriptstyle\frac{{\gamma}_n}{a_n}}I_q$};
\draw[color=black, thick] (4.5,8)node {${\scriptstyle\frac{{\gamma}_1}{a_1}}I_q$};
\draw[color=black, thick] (4.5,1.5)node {${\scriptstyle\frac{1}{{\gamma}_n}}F_R$};
\draw[color=black, thick] (4.5,6.5)node {${\scriptstyle\frac{1}{{\gamma}_1}}F_R$};
\draw[color=black, thick] (4.5,9.2)node {${\scriptstyle d_L}I_q$};
\draw[color=black, thick] (4.0,5.0)node {$\bullet$};
\draw[color=black, thick] (4.0,5.3)node {$\bullet$};
\draw[color=black, thick] (4.0,4.7)node {$\bullet$};
\draw[color=black, thick] (2.55,3.0)node[above] {+};
\draw[color=black, thick] (2.55,8.0)node[above] {+};
\draw[color=black, thick] (3.0,2.60)node[right] {-};
\draw[color=black, thick] (3.0,7.60)node[right] {-};
\draw[color=black, thick] (6.5,5.40)node[right] {+};
\draw[color=black, thick] (6.5,5.40)node[left] {+};
\draw[color=black, thick] (7.6,5)node[above] {+};
\draw[color=black, thick] (5.9,4.9)node[above] {+};
\draw[color=black, thick] (5.9,5)node[below] {+};
\draw[color=black, thick] (6.55,4.50)node[left] {+};
\draw[color=black, thick] (8.0,5.40)node[right] {+};
\end{tikzpicture}
\caption{${\rm Out}=f_L\left(F_R(z)\right)\cdot{\rm In}$\quad
with\quad $f_L(z)=
{\scriptstyle d_L}+\sum\limits_{j=1}^n{\scriptstyle\gamma}_j(z+a_j)^{-1}
$}
\label{Fig:CompositionFeedbackLoop}
\end{figure}

\section{Composition of functions - second version}
\label{Sec:SecondComposition}
\setcounter{equation}{0}

Here, we address a second version of composition of realizations.
Specifically, using \eqref{eq:LeftRightFunctions} we set
\begin{equation}\label{eq:DefSondComposition}
F_L(F_R(z))=D_L+C_L(F_R(z)-A_L)^{-1}B_L.
\end{equation}
To this end, we assume that
\[
n,~ {\rm the~McMillan~degree~of}~ F_L(z), ~{\rm is~
equal~to~the~dimension~of}~ F_R(z).
\]
Moreover assume that the
$n\times n$ matrix $D_R-A_L$ is non-singular, i.e.
\begin{equation}\label{eq:DetCond}
\det(D_R-A_L)\not=0.
\end{equation}
We can now present the main result of this section.

\begin{proposition}\label{Pn:CompositionTensor}
Under the above premises, the composed function, see
\eqref{eq:FrameworkRealization}, \eqref{eq:DefCompFunct},
\eqref{eq:DefSondComposition} and \eqref{eq:DetCond},
is \mbox{$p\times p$-valued} and of McMillan degree $m$.
A corresponding realization array is given by
\[
\begin{footnotesize}
\left(\begin{array}{r|r}
A_{\rm comp}&B_{\rm comp}\\
\hline
C_{\rm comp}&D_{\rm comp}
\end{array}\right)\end{footnotesize}
=
\begin{footnotesize}
\left(\begin{array}{r|r}
A_R-B_R(D_R-A_L)^{-1}C_R&
B_R(D_R-A_L)^{-1}B_L\\
\hline
-C_L(D_R-A_L)^{-1}C_R&
D_L+C_L(D_R-A_L)^{-1}B_L
\end{array}\right).\end{footnotesize}
\]
\end{proposition}

\begin{remark}\label{Rm:RealizeSecondComp}
{\rm
Sometimes we shall find it convenient to re-write the result of
Proposition \ref{Pn:CompositionTensor} as
\[
\left(\begin{smallmatrix}
A_{\rm comp}&B_{\rm comp}\\
C_{\rm comp}&D_{\rm comp}
\end{smallmatrix}\right)
=
\left(\begin{smallmatrix}
A_R&0\\
0&D_L
\end{smallmatrix}\right)
+
\left(\begin{smallmatrix}
B_R\\~\\ C_L
\end{smallmatrix}\right)
\left(\begin{smallmatrix}
D_R-A_L
\end{smallmatrix}\right)^{-1}
\left(\begin{smallmatrix}
-C_R&&B_L
\end{smallmatrix}\right).
\]
}
\end{remark}

\begin{remark}\label{Rm:CondRealizationDependent}
{\rm Note that by \eqref{eq:Coorninates} the matrix $A_L$ is
coordinates-dependent. This in particular implies that
almost always one can make condition \eqref{eq:DetCond} satisfied.
}
\end{remark}

\begin{remark}
{\rm A simple example illustrating the difference between the two
versions of composition dealt with in this work, is when
$F_R(z)=f_R(z)I_n$, where $f_R(z)$ is scalar rational.}
\end{remark}

\noindent
{\bf Proof of Proposition \ref{Pn:CompositionTensor}:} 
In the sequel we shall use the identity,
\begin{equation}\label{eq:Identity}
(I_n+XY)^{-1}=I_n-X(I_m+YX)^{-1}Y \quad\quad\quad
\begin{smallmatrix}
X\in{\mathbb C}^{n\times m}\\~\\
Y\in{\mathbb C}^{m\times n}\\~\\
-1\not\in{\rm spect}(XY).
\end{smallmatrix}
\end{equation}
We now have
\[
\begin{split}
F_L(F_R(z))
&=D_L+C_L\left(\underbrace{C_R(zI_m-A_R)^{-1}B_R
+D_R}_{F_R(z)}-A_L\right)^{-1}B_R\\
&=D_L+C_L\left(C_R(zI_m-A_R)^{-1}B_R
+(D_R-A_L)\right)^{-1}B_R\\
&=D_L+C_L(D_R-A_L)^{-1}\left(I_n+
\underbrace{C_R(zI_m-A_R)^{-1}}_X
\underbrace{B_R(D_R-A_L)^{-1}}_Y\right)^{-1}B_L
\\
&=D_L+C_L(D_R-A_L)^{-1}\times\\
&\hspace{-15mm}\times\left(I_n-
\underbrace{C_R(zI_m-A_R)^{-1}}_X
\left(I_m+
\underbrace{B_R(D_R-A_L)^{-1}}_Y
\underbrace{C_R(zI_m-A_R)^{-1}}_X\right)^{-1}
\underbrace{B_R(D_R-A_L)^{-1}}_Y\right)B_L
\\
&=D_L+C_L(D_R-A_L)^{-1}\times\\
&\hspace{5mm}\times\left(I_n-
C_R\left(zI_m-A_R+
B_R(D_R-A_L)^{-1}C_R\right)^{-1}
B_R(D_R-A_L)^{-1}\right)B_L\\
&=\underbrace{D_L+C_L(D_R-A_L)^{-1}B_L}_{D_{\rm comp}}-\\
&\hspace{5mm}\underbrace{-C_L(D_R-A_L)^{-1}C_R}_{C_{\rm comp}}
\left(zI\underbrace{-A_R+B_R(D_R-A_L)^{-1}C_R}_{-A_{\rm comp}}
\right)^{-1}\underbrace{B_R(D_R-A_L)^{-1}B_L}_{B_{\rm comp}}~,
\end{split}
\]
A critical part is when one substitutes in 
\eqref{eq:Identity}
the values, $~X:=C_R(zI_m-A_R)^{-1}$ and $~Y:=B_R(D_R-A_L)^{-1}$.
\mbox{}\qed\mbox{}\\

We conclude this section by examining the extent to which the main
result is coordinates-dependent.

\begin{remark}{\rm
Assume the realizations of $F_L$ and $F_R$ are minimal and consider a
change of coordinates as in \eqref{eq:Coorninates}, i.e. for some
non-singular $S_L$, $S_R$ ($n\times n$ and $m\times m$ respectively),
\[
\left(\begin{smallmatrix}S_L&0\\0&I_p\end{smallmatrix}\right)^{-1}
\left(\begin{smallmatrix}A_L&B_L\\C_L&D_L\end{smallmatrix}\right)
\left(\begin{smallmatrix}
S_L&0\\0&I_p
\end{smallmatrix}\right)
\quad
\quad
\quad
\left(\begin{smallmatrix}
S_R&0\\0&I_n
\end{smallmatrix}\right)^{-1}
\left(\begin{smallmatrix}
A_R&B_R\\
C_R&D_R
\end{smallmatrix}\right)
\left(\begin{smallmatrix}
S_R&0\\0&I_n
\end{smallmatrix}\right).
\]
Substituting in Proposition \ref{Pn:CompositionTensor} yields
\begin{equation}\label{eq:CoordinateComp}
\begin{smallmatrix}
A_{\rm comp}=&
S_R^{-1}A_RS_R-S_R^{-1}B_R(D_R-S_L^{-1}A_LS_L)^{-1}C_RS_R
&=&
S_R^{-1}\left(A_R-B_R(D_R-S_L^{-1}A_LS_L)^{-1}C_R\right)S_R
\\~\\
B_{\rm comp}=&
S_R^{-1}B_R(D_R-S_L^{-1}A_LS_L)^{-1}S_L^{-1}B_L
&=&
S_R^{-1}\left(B_R(D_R-S_L^{-1}A_LS_L)^{-1}S_L^{-1}B_L\right)
\\~\\
C_{\rm comp}=&
-C_LS_L(D_R-S_L^{-1}A_LS_L)^{-1}C_RS_R
&=&
\left(-C_LS_L(D_R-S_L^{-1}A_LS_L)^{-1}C_R\right)S_R
\\~\\
D_{\rm comp}=&D_L+C_LS_L(D_R-S_L^{-1}A_LS_L)^{-1}S_L^{-1}B_L
&=&
D_L+C_L(S_LD_RS_L^{-1}-A_L)^{-1}B_L~,
\end{smallmatrix}
\end{equation}
which may be a different system.\\
In the special case where,
\begin{equation}
\label{group}
\left(\begin{smallmatrix}S_L&0\\0&I_n\end{smallmatrix}\right)
\left(\begin{smallmatrix}A_L&B_L\\C_L&D_L\end{smallmatrix}\right)
=
\left(\begin{smallmatrix}A_L&B_L\\C_L&D_L\end{smallmatrix}\right)
\left(\begin{smallmatrix}S_L&0\\0&I_n\end{smallmatrix}\right),
\end{equation}
\eqref{eq:CoordinateComp} can be written as the following
change of coordinates,
\[
\left(\begin{smallmatrix}
S_R&0\\0&I_p
\end{smallmatrix}\right)^{-1}
\left(\begin{smallmatrix}
A_{\rm comp}&B_{\rm comp}\\
C_{\rm comp}&D_{\rm comp}
\end{smallmatrix}\right)
\left(\begin{smallmatrix}
S_R&0\\0&I_p
\end{smallmatrix}\right).
\]
We also remark that the set of invertible matrices $S_L$ satisfying
\eqref{group} forms a multiplicative group.
}
\end{remark}

\section{Stieltjes functions}
\label{Sec:Stieltjes}
\setcounter{equation}{0}

Recall that in \eqref{eq:Positive} we described positive functions
$F(z)$ as those that
\begin{equation}\label{eq:PositiveAgain}
\forall z\in\mathbb{C}_r\quad\quad
\begin{matrix}
F(z)\quad{\rm analytic}
\\~\\
\left(F(z)+(F(z))^*\right)
\in\overline{\mathbf P}_p~.
\end{matrix}
\end{equation}
The subset of positive functions in \eqref{eq:PositiveAgain},
where in addition
\begin{equation}\label{eq:DefStieljes}
\forall z\in\mathbb{C}_r\quad\quad
\left(
{\scriptstyle\frac{1}{iz}}F(z)
+
\left(
{\scriptstyle\frac{1}{iz}}F(z)
\right)^*
\right)\in\overline{\mathbf P}_p~,
\end{equation}
are called Stieltjes functions\begin{footnote}{Note that we are
not consistent with \cite[Definition 3.1]{DyukKats1988}
where instead of positive functions described in \eqref{eq:PositiveAgain},
they use Nevanlinna functions analytically mapping the {\em upper}
half plane to itself.}
\end{footnote}.
\vskip 0.2cm

\noindent
In the sequel we shall rely on the following result taken
from\begin{footnote}{A proof of this result is given in
\cite{KreinNudel1977}.}\end{footnote}
\cite[Theorem 3.1]{DyukKats1988} (where originally poles at infinity
are allowed):  Stieltjes functions are exactly those which be can be
written in the form
\begin{equation}\label{eq:integral}
F(z)=i\Delta+\int_0^\infty\frac{z}{t-iz}d\sigma(t),
\quad\Delta\in\mathbf{P}_p\quad\forall z\in\mathbb{C}_r~,
\end{equation}
where the $p\times p$-valued positive measure $\sigma$ satisfies 
\[
\int_0^\infty \frac{d\sigma(t)}{1+t}<\infty.
\]
Here we focus on the rational case, namely where the measure
$\sigma$ has a finite number of jumps.
\vskip 0.2cm

\noindent
For example, a straightforward calculation reveals that all
scalar rational Stieltjes functions of degree one, $f(z)$, 
may be parametrized as,
\[
f(z)=i\left({\scriptstyle\delta+\frac{\beta}{\alpha}}\right)+\frac{
{\scriptstyle\beta}}{z+i{\scriptstyle\alpha}}
\quad\quad
\begin{smallmatrix}
\alpha>0\\~\\
\beta>0\\~\\
\delta\geq 0.
\end{smallmatrix}
\]
This observation is  next generalized to all {\em rational} functions
satisfying \eqref{eq:integral}. 

\begin{proposition}\label{Pn:RealizationStieltjes}
Let $F(z)$ be a $p\times p$-valued rational function, analytic
at the origin and at infinity, of McMillan degree $n$.\\
$F(z)$ is a Stieltjes function, satisfying
\eqref{eq:PositiveAgain} and \eqref{eq:DefStieljes},
if and only if, it can be written as,
\begin{equation}\label{eq:RationalStieljes}
F(z)=i\left(
C{\alpha}^{-1}C^*+\delta
\right)+
C(zI_n+i\alpha)^{-1}C^*
\quad\quad{\rm with}\quad\quad
\begin{smallmatrix}
C\in\mathbb{C}^{p\times n}~{\rm full~rank}
\\~\\
\alpha\in{\mathbf P}_n
\\~\\
\delta
\in\overline{\mathbf P}_p~.
\end{smallmatrix}
\end{equation}
\end{proposition}

\noindent
{\bf Proof :}\quad First recall, see e.g. \cite[Lemma 1.1(II)]{AlpayLew2011},
\cite[Chapter 5]{AnderVongpa1973}, \cite{DelGeinKamp1981}, that from
the realization matrix formulation of the K-Y-P Lemma it follows that
a rational function $F(z)$, analytic at infinity, is positive\begin{footnote}
{It may be complex or real.}\end{footnote} if and only if, up to change of
coordinates, its minimal realization satisfies
\begin{equation}\label{eq:KYPgeneral}
\left(\begin{smallmatrix}-I_n&0\\0&I_p\end{smallmatrix}\right)
R_F
+
R_F^*
\left(\begin{smallmatrix}-I_n&0\\0&I_p\end{smallmatrix}\right)
\in\overline{\mathbf P}_{n+p}~.
\end{equation}
Next, note that \eqref{eq:integral} gives an analytic extension of
$F(z)$
to
$\mathbb{C}\smallsetminus i\mathbb R_-$ such that
\begin{equation}\label{asdf}
(F
(-
{z}^*))^*=-F
(z).
\end{equation}
Recall now that positive functions which in addition satisfy
\eqref{asdf} are called Positive Odd. In the real rational
case they are known in electrical engineering as {\em Lossless}
or {\em Foster}, see e.g. \cite{AnderVongpa1973},
\cite{Belev1968}, \cite{Wohl1969}.
\vskip 0.2cm

\noindent
Furthermore, if in addition $F(z)$ is odd,
i.e.
\eqref{asdf} holds, then its realization array $R_F$ may be chosen so that
\begin{equation}\label{eq:KYPforPRO}
\left(\begin{smallmatrix}-I_n&0\\0&I_p\end{smallmatrix}\right)
R_F+R_F^*
\left(\begin{smallmatrix}-I_n&0\\0&I_p\end{smallmatrix}\right)
=0,
\end{equation}
see e.g. \cite[Theorem 4.1]{AlpayGohberg1988},
\cite[Section 5.2]{AnderVongpa1973}. Note now that
\eqref{eq:KYPforPRO} means that the \mbox{$(n+p)\times(n+p)$} matrix 
$\left(\begin{smallmatrix}-I_n&0\\0&I_p\end{smallmatrix}\right)
R_f$ is skew-Hermitian, namely,
\[
\left(
i
\left(\begin{smallmatrix}-I_n&0\\0&I_p\end{smallmatrix}\right)
R_F
\right)
=
\left(
i
\left(\begin{smallmatrix}-I_n&0\\0&I_p\end{smallmatrix}\right)
R_F
\right)^*,
\]
which in turn can be written as,
\[
R_F=\begin{footnotesize}\left(\begin{array}{c|c}
-i\alpha&C^*\\
\hline
C&i\Delta\end{array}\right)
\end{footnotesize}
\quad\quad{\rm with}\quad\quad
\begin{smallmatrix}
\alpha={\alpha}^*\\~\\
C\in\mathbb{C}^{p\times n}\\~\\
\Delta={\Delta}^*.
\end{smallmatrix}
\]
Thus far one can conclude that
\begin{equation}\label{eq:BasicRealization1}
F(z)=i\Delta+C(zI_n+i\alpha)^{-1}C^*
\quad\quad{\rm with}\quad\quad
\begin{smallmatrix}
\alpha={\alpha}^*\\~\\
C\in\mathbb{C}^{p\times n}\\~\\
\Delta={\Delta}^*.
\end{smallmatrix}
\end{equation}
We next show that
\begin{equation}\label{eq:CondAlpha}
\alpha\in\mathbf{P}_n~,
\end{equation}
and that
\begin{equation}\label{eq:AlternativeCondAlpha}
{\scriptstyle\Delta}=C{\alpha}^{-1}C^*+{\scriptstyle\delta}
\quad\quad{\rm for~some}\quad\quad
{\scriptstyle\delta}\in\overline{\mathbf P}_p~.
\end{equation}
To this end, using the fact that by assumption, $F(z)$ is analytic at
the origin, i.e. $\alpha$ is non-singular,
we shall find it convenient to re-write the $F(z)$ in
\eqref{eq:BasicRealization1} as
\[
F(z)
=i\left(\Delta-C{\alpha}^{-1}C^*\right)+izC
{\alpha}^{-1}(z{\alpha}^{-1}+iI_n)^{-1}{\alpha}^{-1}C^*
\quad\quad{\rm with}\quad\quad
\begin{smallmatrix}
\alpha={\alpha}^*\\~\\
C\in\mathbb{C}^{p\times n}\\~\\
\Delta={\Delta}^*,
\end{smallmatrix}
\]
and hence,
\[
{\scriptstyle\frac{1}{iz}}F(z)=
{\scriptstyle\frac{1}{z}}(\Delta-C{\alpha}^{-1}C^*)
+C{\alpha}^{-1}(z{\alpha}^{-1}+iI_n)^{-1}{\alpha}^{-1}C^*.
\]
We can now substitute the above $F(z)$ in \eqref{eq:DefStieljes}
to obtain,
\[
{\scriptstyle
\left(
{\scriptstyle\frac{1}{z}}(\Delta-C{\alpha}^{-1}C^*)
+C{\alpha}^{-1}(z{\alpha}^{-1}+iI_n)^{-1}{\alpha}^{-1}C^*
+
\left(
{\scriptstyle\frac{1}{z}}(\Delta-C{\alpha}^{-1}C^*)
+C{\alpha}^{-1}(z{\alpha}^{-1}+iI_n)^{-1}{\alpha}^{-1}C^*
\right)^*
\right)}\in\overline{\mathbf P}_p
\quad\quad\forall z\in\mathbb{C}_r~,
\]
i.e.
\begin{equation}\label{eq:XX}
{\scriptstyle
2{\rm Re}(z)
\left(
{\scriptstyle\frac{1}{|z|^2}}(
\Delta
-C{\alpha}^{-1}C^*)
+
\left(C{\alpha}^{-1}(z{\alpha}^{-1}+iI_n)^{-1}\right)
{\alpha}^{-1}
\left(C{\alpha}^{-1}(z{\alpha}^{-1}+iI_n)^{-1}\right)^*
\right)
}\in\overline{\mathbf P}_p
\quad\quad\quad
\forall z\in\mathbb{C}_r~.
\end{equation}
Clearly, having \eqref{eq:CondAlpha} along with
\eqref{eq:AlternativeCondAlpha} implies that \eqref{eq:XX} holds. Thus,
there is the converse direction to consider.
\vskip 0.2cm

\noindent
First, note that since  \eqref{eq:XX} holds in particular
for all points of $z\in\mathbb{C}_r$ (up to $n$ points) so that the
matrix $z{\alpha}^{-1}+iI_n$ is nearly singular, this in fact implies that
${\alpha}^{-1}\in\mathbf{P}_n$ i.e. \eqref{eq:CondAlpha} holds.
\vskip 0.2cm

\noindent
Similarly, as \eqref{eq:XX} is satisfied in particular for $z\in\mathbb{C}_r$
``sufficiently small", it implies that \eqref{eq:AlternativeCondAlpha}
holds as well, so the claim is established.
\qed
\vskip 0.2cm

\begin{remark}\label{Rm:FirstPrarametrizationStieljes}
Eq. \eqref{eq:RationalStieljes} may be viewed as a
parametrization of all rational Stieltjes function
analytic at the origin and at infinity.
\end{remark}
\vskip 0.2cm

\noindent
We now next review this result. To this end we recall the following.

\begin{remark}\label{Rm:PosSemiDef}
Consider the following statements
for a full-rank matrix ${\scriptstyle Z}\in{\mathbb C}^{p\times n}$.
\vskip 0.2cm

\begin{itemize}
\item[(i)~~~]{}\quad$
\left(\begin{smallmatrix}
X&Z^*\\Z&Y
\end{smallmatrix}\right)\in\overline{\mathbf P}_{n+p}~.
$
\vskip 0.2cm

\item[(ii)~~]{}\quad$
{\scriptstyle Y}\in{\mathbf P}_p\quad{\rm and}\quad
{\scriptstyle X}-{\scriptstyle Z^*Y^{-1}Z}\in\overline{\mathbf P}_n~.$
\vskip 0.2cm

\item[(iii)~]{}\quad$
{\scriptstyle X}\in{\mathbf P}_n\quad{\rm and}\quad
{\scriptstyle Y}-{\scriptstyle ZX^{-1}Z^*}\in\overline{\mathbf P}_p~.$
\end{itemize}
\vskip 0.2cm

\noindent
Then, (ii) implies (i) and if $n\geq p$ then the converse is true as well.
\vskip 0.2cm

\noindent
Then, (iii) implies (i) and if $p\geq n$ then the converse is true as well.
\end{remark}

Remark \ref{Rm:PosSemiDef} leads to the conclusion that
Proposition \ref{Pn:RealizationStieltjes} and Remark
\ref{Rm:FirstPrarametrizationStieljes} can be alternatively formulated
as follows.

\begin{remark}\label{Rm:AlternativeCharacterizationStieljes}
All $n\times n$-valued rational Stieltjes function $F(z)$, analytic
at the origin and at infinity, of McMMillan degree $m$, with
$n\geq m$, may be parametrized as
\[
F(z)=i\Delta
+C\left(zI_m+i(C^*{\Delta}^{-1}C+{\scriptstyle\eta})\right)^{-1}C^*
\quad\quad\quad
\begin{smallmatrix}
\Delta\in{\mathbf P}_n
\\~\\
C\in\mathbb{C}^{n\times m}~{\rm full~rank}
\\~\\
\eta\in\overline{\mathbf P}_m~,
\end{smallmatrix}
\]
namely admitting a realization of the form,
\[
R_F=\begin{footnotesize}\left(\begin{array}{c|c}
-i(C^*{\Delta}^{-1}C+{\scriptstyle\eta})&C^*\\
\hline
C&i\Delta
\end{array}\right)\end{footnotesize}
\quad\quad\quad
\begin{smallmatrix}
\Delta\in{\mathbf P}_n
\\~\\
C\in\mathbb{C}^{n\times m}~{\rm full~rank}
\\~\\
\eta\in\overline{\mathbf P}_m~.
\end{smallmatrix}
\]
\end{remark}

In the sequel, we shall find it convenient to use the following.

\begin{remark}\label{Rm:AlternativeRealizationStieljes}
Denoting
\[
\gamma:=-iC,
\]
one can re-write the realization
of $F(z)$ in Proposition \ref{Pn:RealizationStieltjes} as
\begin{equation}\label{eq:RealizationArrayStieljes1}
R_F=
\begin{footnotesize}\left(\begin{array}{c|c}
-i\alpha&C^*\\
\hline
C&i\Delta
\end{array}\right)\end{footnotesize}
=
i
\left(
\begin{smallmatrix}
-I_n&0\\0&I_p
\end{smallmatrix}
\right)
\left(
\begin{smallmatrix}
\alpha&{\gamma}^*\\
\gamma&\Delta
\end{smallmatrix}
\right)
\end{equation}
where
\begin{equation}\label{eq:RealizationArrayStieljes2}
\left(
\begin{smallmatrix}
\alpha&{\gamma}^*\\
\gamma&\Delta
\end{smallmatrix}
\right)
\in\overline{\mathbf P}_{n+p}
\quad\quad{\rm and}\quad\quad
\begin{smallmatrix}
\alpha\in\mathbf{P}_n&&~\\~\\
\gamma\in\mathbb{C}^{p\times n}&&{\rm full~rank}\\~\\
n\geq p.
\end{smallmatrix}
\end{equation}
\end{remark}
\vskip 0.2cm

\noindent
To further emphasize the difference between Stieltjes functions
and those discussed in
Section \ref{Sec:Motivating applications}, we have the following.

\begin{remark}\label{Rm:Invertible}
As already mentioned the family of Stieltjes functions is a
proper subset of Positive Odd functions\begin{footnote}{
For example $\frac{1}{s+1}$ is a positive function which is not
odd and
$\frac{1}{s+i}$ is a positive odd function which is not
Stieltjes.}\end{footnote}, which in turn is a
proper subset of Positive functions. Consider the following
properties.

\noindent
$\bullet$\quad
Each of these three sets is closed under positive scaling and summation
and thus is a convex cone.
\vskip 0.2cm

$\bullet$\quad
Both Positive functions and its subset of Positive Odd functions
are closed under inversion, namely if $F(z)$ is Positive (Odd) then
$\left(F(z)\right)^{-1}$ is well defined and is Positive (Odd). Thus,
each of these two sets is a Convex Invertible Cone. In
\cite{CohenLew2007} this fact was explored in the framework of real
functions.
\vskip 0.2cm

$\bullet$\quad
If $F(z)$ is a Stieltjes function then $\left(F(z)\right)^{-1}$
is well defined Positive Odd function, which can not be a Stieltjes
function. Indeed, from Remark \ref{Rm:AlternativeRealizationStieljes}
it follows that
\[
-i\left(\lim\limits_{z~\rightarrow~\infty}F(z)\right)\in{\mathbf P}_p
\quad\quad{\rm but}\quad\quad+i\left(
\lim\limits_{z~\rightarrow~\infty}
\left(F(z)\right)^{-1}\right)\in{\mathbf P}_p~.
\]
\end{remark}
\vskip 0.2cm

\noindent
Recall now that the K-Y-P Lemma, see e.g. \cite{AlpayLew2011}
\cite[Chapter 5]{AnderVongpa1973}, characterizes positive 
rational functions, along with some sub-families, through
properties of their minimal realizations. We can now introduce
an adaptation of the K-Y-P Lemma to Stieltjes functions and
then use it to construct from a given realization ~{\em  a whole family}
of Stieltjes functions of various dimensions and McMillan degrees.

\begin{corollary}\label{Cy:AlternativeDescriptionStieljes}
Let $F(z)$ be a rational function with $p$ outputs, analytic at infinity,
of McMillan degree $n$, as in
\eqref{eq:BasicRealization} and \eqref{eq:RealizationArray} i.e.
\[
F(z)=D+C(zI_n-A)^{-1}B\quad\quad\quad
R_F=
\begin{footnotesize}
\left(\begin{array}{c|c}
A&B\\
\hline
C&D\end{array}\right).
\end{footnotesize}
\]
Then, $F(z)$ is a Stieltjes function, if and only if the realization
$R_F$ can be chosen so that each of the four blocks $A, B, C, D$ is
of a full rank and
\begin{equation}\label{eq:CharacterizeStieljes}
-i\left(\begin{smallmatrix}-I_n&0\\~~0&I_p
\end{smallmatrix}\right)R_F\in\overline{\mathbf P}_{n+p}~.
\end{equation}
Moreover, let $U\in\mathbb{C}^{\nu\times n}$ and
$V\in\mathbb{C}^{\pi\times p}$ be full rank matrices,
for some $\nu\in[1, n]$ and $\pi\in[1, p]$.\quad
Then,
\[
R_{\hat{F}}=
\left(\begin{smallmatrix}U&0\\0&V\end{smallmatrix}\right)
R_F
\left(\begin{smallmatrix}U&0\\0&V\end{smallmatrix}\right)^*,
\]
is a realization of a \mbox{$\pi\times\pi$-valued} Stieltjes function 
$\hat{F}(z)$ (analytic at the origin) of McMillan degree $\nu$.
\end{corollary}

\noindent
Indeed the claim follows from Proposition \ref{Pn:RealizationStieltjes},
Remarks \ref{Rm:PosSemiDef} and \ref{Rm:AlternativeRealizationStieljes}
along with the fact that $T$ is a full rank matrix so that the $H=TMT^*$
is well defined, the product matrix $H$ is positive (semi-)definite if
and only if $M$ is positive (semi-)definite.
\vskip 0.2cm

\noindent
Following Remark \ref{Rm:AlternativeRealizationStieljes} and
Corollary \ref{Cy:AlternativeDescriptionStieljes} we shall
call a realization $R_F$ of a Stieltjes function ~{\em canonical} if 
it satisfies \eqref{eq:RealizationArrayStieljes1} with
\eqref{eq:RealizationArrayStieljes2}, or equivalently
\eqref{eq:CharacterizeStieljes}.
\vskip 0.2cm

\noindent
Note that a realization remains {\em canonical} under {\em unitary}
change of coordinates, i.e. in \eqref{eq:Coorninates} $S^{-1}=S^*$.
\vskip 0.2cm

\noindent
We next show that the set of rational Stieltjes function is closed under
the second version of composition of functions, see Section
\ref{Sec:SecondComposition}.

\begin{proposition}\label{Pn:CompStieljes}
Consider a pair of rational Stieltjes functions, $F_L(z)$, $F_R(z)$,
analytic at the origin and at infinity: $F_L(z)$ is
\mbox{$p\times p$-valued} of McMillan degree $n$ and $F_R(z)$ is
\mbox{$n\times n$-valued} of McMillan degree $m$ with $n\geq m$.
\vskip 0.2cm

\noindent
If the realization of both $F_L(z)$ and $F_R(z)$ is
canonical\begin{footnote}{Thus the condition in \eqref{eq:DetCond} is
trivially satisfied.}\end{footnote},
then the composition $F_L(F_R)$ in \eqref{eq:DefSondComposition} is a
$p\times p$-valued rational Stieltjes function of McMillan degree $m$,
given in a canonical realization.
\end{proposition}
\vskip 0.4cm

\noindent
{\bf Proof :}\quad
Using Proposition \ref{Pn:RealizationStieltjes} along with Remark
\ref{Rm:AlternativeRealizationStieljes} below left, we have $R_{F_L}$
a realization array of $F_L(z)$ and using Remark
\ref{Rm:AlternativeCharacterizationStieljes}, below right, we have
$R_{F_R}$ a realization array of $F_R(z)$,
\[
\begin{smallmatrix}
R_{F_L}=
\begin{footnotesize}
\left(\begin{array}{c|c}
-i{\alpha}_L&(i{\gamma}_L)^*\\
\hline
i{\gamma}_L&i(
{\gamma}_L
{{\alpha}_L}^{-1}
{{\gamma}_L}^*
+\delta)
\end{array}\right)
\end{footnotesize}
&&&&
R_{F_R}=
\begin{footnotesize}
\left(\begin{array}{c|c}
-i(
{{\gamma}_R}^*
{{\scriptstyle\Delta}_R}^{-1}
{\gamma}_R
+\eta
)&(i{\gamma}_R)^*\\
\hline
i{\gamma}_R&i{\Delta}_R\end{array}\right)
\end{footnotesize}
\\~\\~\\
{\scriptstyle{\gamma}_L}\in\mathbb{C}^{p\times n}
\quad{\rm full~rank}
&&&&
{\scriptstyle{\gamma}_R}\in\mathbb{C}^{n\times m}
\quad{\rm full~rank}
\\~\\
{\scriptstyle{\alpha}_L}\in\mathbf{P}_n
&&&&
{{\scriptstyle\Delta}_R}\in\mathbf{P}_n
\\~\\
\delta\in\overline{\mathbf P}_p
&&&&
\eta\in\overline{\mathbf P}_m~.
\end{smallmatrix}
\]
Substituting in Proposition \ref{Pn:CompositionTensor} yields that a
realization of $F_{\rm comp}=F_L(F_R)$ in
\eqref{eq:DefSondComposition} is given by
\[
\begin{matrix}
R_{F_{\rm comp}}=
{\scriptstyle i}
\begin{footnotesize}\left(\begin{array}{c|c}
-({\gamma}_R^*{{\scriptstyle\Delta}_R}^{-1}{\gamma}_R+\eta)
+{\gamma}_R^*
\left(\begin{smallmatrix}{\Delta}_R+{\alpha}_L\end{smallmatrix}\right)^{-1}
{\gamma}_R
&
{\gamma}_R^*
\left(\begin{smallmatrix}{\Delta}_R+{\alpha}_L\end{smallmatrix}\right)^{-1}
{\gamma_L}^*\\
\hline
-{\gamma}_L
\left(\begin{smallmatrix}{\Delta}_R+{\alpha}_L\end{smallmatrix}\right)^{-1}
{\gamma}_R
&
{\gamma}_L
{{\alpha}_L}^{-1}
{{\gamma}_L}^*
+\delta
-{\gamma}_L
\left(\begin{smallmatrix}{\Delta}_R+{\alpha}_L\end{smallmatrix}\right)^{-1}
{\gamma_L}^*
\end{array}\right)
\end{footnotesize}
\end{matrix}
\]
A straightforward exercise enables one to re-write this realization as,
\begin{equation}
R_{F_{\rm comp}}={\scriptstyle i}
\left(\begin{smallmatrix}-I_m&0\\0&I_p\end{smallmatrix}\right)
\underbrace{
\left(
\left(\begin{smallmatrix}\eta&0\\0&\delta\end{smallmatrix}\right)
+
\left(\begin{smallmatrix}{\gamma}_R^*\\
-{\gamma}_L{\alpha}_L^{-1}{\Delta}_R\end{smallmatrix}\right)
\underbrace{
\left({\scriptstyle{\Delta}_R}^{-1}-\left(\begin{smallmatrix}
{\Delta}_R+{\alpha}_L\end{smallmatrix}\right)^{-1}\right)
}_M
\left(\begin{smallmatrix}{\gamma}_R^*\\
-{\gamma}_L{\alpha}_L^{-1}{\Delta}_R\end{smallmatrix}\right)^*
\right).
}_W
\end{equation}
Now as by assumption, both ${\alpha}_L$ and ${\scriptstyle{\Delta}_R}$
are in ${\mathbf P}_n$, it implies that
\[
{\scriptstyle M}:=\left({\scriptstyle{\Delta}_R}^{-1}-
\left(\begin{smallmatrix}
{\Delta}_R+{\alpha}_L\end{smallmatrix}\right)^{-1}\right)
\in{\mathbf P}_n~,
\]
as well. From the structure it follows that\begin{footnote}{To be
precise, ${\rm rank}(W)=\min(n,~m+p)$.}\end{footnote}
\[
{\scriptstyle W}:=
\left(
\left(\begin{smallmatrix}\eta&0\\0&\delta\end{smallmatrix}\right)
+
\left(\begin{smallmatrix}{\gamma}_R^*\\
-{\gamma}_L{\alpha}_L^{-1}{\Delta}_R\end{smallmatrix}\right)
{\scriptstyle M}
\left(\begin{smallmatrix}{\gamma}_R^*\\
-{\gamma}_L{\alpha}_L^{-1}{\Delta}_R\end{smallmatrix}\right)^*
\right)\in\overline{\mathbf P}_{m+p}~.
\]
From Remarks \ref{Rm:PosSemiDef} and
\ref{Rm:AlternativeRealizationStieljes} it thus follows that the
resulting \mbox{$F_{\rm comp}=F_L(F_R)$} in
\eqref{eq:DefSondComposition}, is a Stieltjes function.
\qed

\section{Future work}
\label{Sec:ConclRem}
\setcounter{equation}{0}

In this work we focused on composition of rational functions
their state-space realization and applications.
Yet, this research area
is mostly open. We here
point out at three sample problems of various level of importance.
\vskip 0.2cm

\noindent
$\bullet\quad$
Assume having a small set of ``simple" (e.g. low degree)
rational functions as ``building blocks".\\
Synthesis: What functions can be generated from these
building blocks.\\
Analysis: Given a complicated rational function, can it be, and if
yes how, ``factorized" or ``decomposed" into a composition of
simpler building blocks.
\vskip 0.2cm

\noindent
Synthesis is further discussed below. To emphasize
the the importance of analysis, recall that in Remark \ref{Rm:degree}
it was pointed out that in the first version composition,
the McMillan degree of $F_L(F_R)$ is equal to the ~{\em product}
of the McMillan degrees of the original functions
$F_L(z)$ and $F_R(z)$. Thus ``decomposition" may significantly simplify
the functions at hand.
\vskip 0.2cm

\noindent
As an illustration consider the following rational function of two
variables
\begin{equation}\label{eq:Phi}
\phi(F, G):=\left(F+G^{-1}\right)^{-1}.
\end{equation}
Note that this is function on {\em non-commuting} variables, in fact,
\[
\phi(G^{-1}, F^{-1})=\phi(F, G).
\]
Note now that the driving point impedance in Figure
\ref{Fig:CircuitFeedbackLoop4} can be written as,
\[
{\rm Z}_{\rm in}=\phi(Y_F,~Z_G),
\]
and a basic feedback loop in Figure \ref{Figure:FeedbackLoop1},
may be viewed as
\[
{\rm Out}=\phi(F, G)\cdot{\rm In}.
\]
Now composition of such $\phi$  functions yields
\[
\phi(F_c, G_c)={\scriptstyle\left(\underbrace{F_a+G_a^{-1}}_{F_c}
+\underbrace{\left(F_b+G_b^{-1}\right)^{-1} }_{G_c^{-1}}\right)^{-1}}
\quad{\rm where}\quad
\begin{smallmatrix}F_c&=&\phi(F_a, G_a)\\~\\ G_c&=&\phi(F_b, G_b).
\end{smallmatrix}
\]
An illustration of the converse problem let the starting point
be the above function
\[
\left(F_a+G_a^{-1}+\left(F_b+G_b^{-1}\right)^{-1}\right)^{-1},
\]
and using $\phi$ from \eqref{eq:Phi},
one needs to rewrite it in the form of
\[
\phi\left(~\phi(F_a,G_a),~\left(\phi(F_a,G_a)\right)^{-1}\right).
\]
$\bullet\quad$
In Section \ref{Sec:Motivating applications} we presented inter-relations
between (i) positive real rational functions (ii) driving point
inpedance of $R-L-C$ networks and (iii) feedback loops.
\vskip 0.2cm

\noindent
As already mentioned, identifying items (i) with (ii) is classical.
Moreover, there is a whole list of synthesis schemes how to
construct an $R-L-C$ circuit whose driving point impedance realizes
a prescribed positive real rational function: Bott-Duffin, Brune,
Darlington, Foster to name but few, see e.g. \cite{AnderVongpa1973},
\cite{Chen1964}.
\vskip 0.2cm

\noindent
The inter-relation between (i) and (iii) suggests that one can
exploit these electrical circuit synthesis schemes to construct, out of
simple building blocks, an elaborate network of feedback loops. As
a potential application, see Figures \ref{Fig:CompositionFeedbackLoop}
above or \ref{Fig:MultipleFeedbackLoop} below.

\begin{figure}[ht!]
\begin{tikzpicture}{scale=4.0}
\draw[color=black, thick] [->] (1.5,0.5) -- (3,0.5){};
\draw[color=black, thick] [->] (4,0.5) -- (5.5,0.5){};
\draw[color=black, thick] [->] (0,3.5) node[left]{In} -- (1.15,3.5) {};
\draw[color=black, thick] [->] (1.85,3.5) -- (3,3.5){};
\draw[color=black, thick] [->] (4,3.5) -- (5.5,3.5){};
\draw[color=black, thick] [->] (5.5,3.5) -- (6.5,3.5) node[right] {Out};
\draw[color=black, thick] [->] (5.5,0.5) -- (5.5,1.65){};
\draw[color=black, thick] [->] (1.5,2) -- (1.5,0.5){};
\draw[color=black, thick] [->]  (1.5,2) -- (1.5,3.15){};
\draw[color=black, thick] [->]  (1.5,5) -- (1.5,3.85){};
\draw[color=black, thick] [->]  (3,5) -- (1.5,5){};
\draw[color=black, thick] [->]  (5.5,5) -- (4,5){};
\draw[color=black, thick] [->]  (5.5,3.5) -- (5.5,5){};
\draw[color=black, thick] [->]  (5.5,3.5) -- (5.5,2.35){};
\draw[color=black, thick] [->]  (3,2) -- (1.5,2){};
\draw[color=black, thick] [->]  (5.15,2) -- (4,2){};
\draw[color=black, thick] [-]  (3,5.35) -- (4,5.35){};
\draw[color=black, thick] [-]  (3,5.35) -- (3,4.65){};
\draw[color=black, thick] [-]  (4,5.35) -- (4,4.65){};
\draw[color=black, thick] [-]  (4,5.35) -- (4,4.65){};
\draw[color=black, thick] [-]  (3,4.65) -- (4,4.65){};
\draw[color=black, thick] [-]  (3,3.85) -- (4,3.85){};
\draw[color=black, thick] [-]  (3,3.15) -- (4,3.15){};
\draw[color=black, thick] [-]  (3,3.15) -- (3,3.85){};
\draw[color=black, thick] [-]  (4,3.15) -- (4,3.85){};
\draw[color=black, thick] [-]  (3,2.35) -- (4,2.35){};
\draw[color=black, thick] [-]  (3,2.35) -- (3,1.65){};
\draw[color=black, thick] [-]  (4,2.35) -- (4,1.65){};
\draw[color=black, thick] [-]  (3,1.65) -- (4,1.65){};
\draw[color=black, thick] [-]  (3,0.85) -- (4,0.85){};
\draw[color=black, thick] [-]  (3,0.15) -- (4,0.15){};
\draw[color=black, thick] [-]  (3,0.15) -- (3,0.85){};
\draw[color=black, thick] [-]  (4,0.15) -- (4,0.85){};
\draw[color=black, thick] (1.5,3.5) circle (0.35){};
\draw[color=black, thick] (5.5,2.0) circle (0.35){};
\draw[color=black, thick] (3.5,5)node {$F_a$};
\draw[color=black, thick] (3.5,3.5)node {$G_a$};
\draw[color=black, thick] (3.5,2)node {$G_b$};
\draw[color=black, thick] (3.5,0.5)node {$F_b$};
\draw[color=black, thick] (1.05,3.5)node[above] {+};
\draw[color=black, thick] (1.5,3.95)node[right] {-};
\draw[color=black, thick] (1.5,3.05)node[right] {-};
\draw[color=black, thick] (5.5,2.45)node[right] {+};
\draw[color=black, thick] (5.5,1.55)node[right] {-};
\end{tikzpicture}
\caption{$
{\rm Out}=\left(F_a+{G_a}^{-1}+\left(F_b+{G_b}^{-1}
\right)^{-1}\right)^{-1}\cdot{\rm In}
$}
\label{Fig:MultipleFeedbackLoop}
\end{figure}
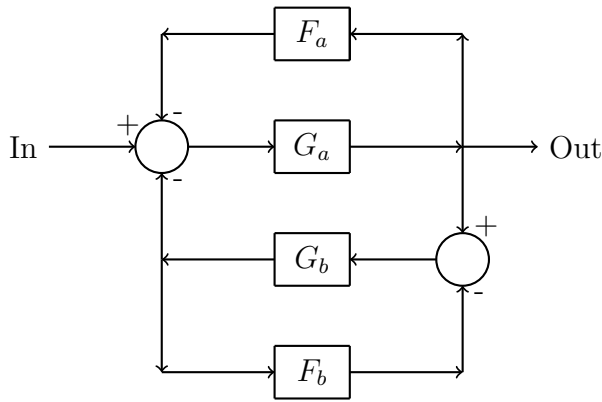
\vskip 0.2cm

\noindent
For example the above $\phi$ in \eqref{eq:Phi} is positive real in the
sense that if $F(z)$ and $G(z)$ are positive real, $\phi$ satisfies
\eqref{eq:Positive} or \eqref{eq:PositiveAgain}. However, this study
requires caution in at least two ways:\\
(i) The application to constructing feedback loop networks, 
{\em transcends} the framework where the building blocks, like $F(z)$
or $G(z)$, are rational positive. For instance, in Figure
\ref{Figure:FeedbackLoop1} the functions $F(z)$ and $G(z)$ are only
required to satisfy ${\rm det}(G)\not\equiv 0$ and
{\mbox{${\rm det}(F+G^{-1})\not\equiv 0$}.
They need not be positive and in principle even not necessarily 
rational.\\
(ii) Study of rational functions of {\em non-commuting} variables
in general and
positive real is praticular, has been flourishing recently, as sample
references see e.g. \cite{AbuAlpColLewSab2018}, \cite{BallGroTerH2018a},
\cite{BallGroTerH2018b}, \cite{BallKal2015}, \cite{LugNed2017},
\cite{PascPassTull2018}, \cite{Popes2013}, \cite{{KaluVin2014}}.
Nevertheless many properties of these functions are yet to be explored. For
example, there is a long way to go to extend (as proposed above) some of the
known electrical circuits synthesis schemes to the framework of
non-commuting variables in order to render it an engineering tool for
designing multi-inputs multi-outouts feedback networks.\\
\vskip 0.2cm

\noindent
$\bullet\quad$
Assuming the dimensions of all matrices involved are suitable and that
$M$ is non-singular, the results in Remarks \ref{Rm:RealizeInverse},
\ref{Rm:RealizeFirstComp} and \ref{Rm:RealizeSecondComp} are all in
the framework of,
\[
\left(\begin{smallmatrix}
Y&0\\0&Z
\end{smallmatrix}\right)
+
\left(\begin{smallmatrix}
\pm{B}_R\\~~C_L
\end{smallmatrix}\right)
{\scriptstyle M}^{-1}
\left(\begin{smallmatrix}
\pm{C}_R&&B_L
\end{smallmatrix}\right).
\]
This observation calls for further investigation.

\end{document}